\documentclass[12pt,reqno]{amsart}
\usepackage{a4wide,euscript,amscd,amssymb,slashed}
\newcommand{\eus}{\EuScript}

\newtheorem{thrm}[equation]{Theorem}
\newtheorem{cor}[equation]{Corollary}
\newtheorem{prop}[equation]{Proposition}
\newtheorem{lemma}[equation]{Lemma}

\theoremstyle{remark}
\newtheorem*{remark}{Remark}
\newtheorem*{remarks}{Remarks}

\newtheoremstyle{citing}{3pt}{3pt}{\itshape}{}{\bfseries}{.}{.5em}{\thmnote{#3}}
\theoremstyle{citing}
\newtheorem*{varthrm}{}

\DeclareMathOperator{\const}{const}
\DeclareMathOperator{\dist}{dist}
\DeclareMathOperator{\id}{id}
\DeclareMathOperator{\re}{Re}
\DeclareMathOperator{\im}{Im}
\DeclareMathOperator{\vol}{vol}

\DeclareMathOperator{\End}{End}
\DeclareMathOperator{\Ker}{Ker}
\DeclareMathOperator{\Coker}{Coker}

\DeclareMathOperator{\Ind}{index}

\newcommand{\CX}{\mathbb{C}}
\newcommand{\CP}{\mathbb{C}P}
\newcommand{\RE}{\mathbb{R}}
\newcommand{\RP}{\mathbb{R}P}
\newcommand{\ZE}{\mathbb{Z}}
\newcommand{\HH}{\mathcal{H}}

\newcommand{\corner}{\lrcorner}
\newcommand{\eps}{\varepsilon}
\newcommand{\minus}{\smallsetminus}
\newcommand{\ga}{g_{a^2}}
\newcommand{\gcyl}{g_{\mathrm{cyl}}}
\newcommand{\Mloc}{M_{\mathrm{loc}}}
\newcommand{\Wloc}{W_{\mathrm{loc}}}
\newcommand{\p}{\partial}

\newcommand{\we}{\wedge}
\newcommand{\pol}{{\textstyle\frac12}}
\newcommand{\phiCY}{\varphi_{CY}}
\newcommand{\phiD}{\varphi_D}
\newcommand{\tV}{\tilde{V}}
\newcommand{\tD}{\tilde{D}}
\newcommand{\vea}{v_{\mathbf{e},L}}
\newcommand{\xa}{X_{a^2}}
\newcommand{\xat}{X_{a^2,\theta}}

\let\phi=\varphi
\let\epsilon=\varepsilon

\begin{document}

\pagestyle{headings}

\title[Coassociative K3 fibrations]{Coassociative K3 fibrations of
    compact $G_2$-manifolds}
\author[A. Kovalev]{Alexei Kovalev}

\address{DPMMS, University of Cambridge,
    Centre for \mbox{Mathematical} Sciences,
    \mbox{Wilberforce} Road, Cambridge CB3 0WB,
    Great Britain}
\email{a.g.kovalev@dpmms.cam.ac.uk}

\begin{abstract}
A class of examples of Riemannian metrics with holonomy $G_2$ on compact
7-manifolds was constructed by the author in~\cite{g2} using a certain
`generalized connected sum' of two asymptotically cylindrical 
manifolds with holonomy $SU(3)$. We consider, on each of the
two initial $SU(3)$-manifolds, a fibration arising from a Lefschetz pencil
of K3 surfaces. The gluing of two such K3 fibrations yields a coassociative
fibration of the connected sum $G_2$-manifold over a 3-dimensional sphere.
The singular fibres of this fibration are diffeomorphic to K3 orbifolds
with ordinary double points and are parameterized by a Hopf-type link
in~$S^3$. We believe that these are the first examples of fibrations of
compact manifolds of holonomy $G_2$ by coassociative minimal submanifolds.
\end{abstract}

\maketitle

\section*{Introduction}

It is well-known that if a K3 surface contains an elliptic curve then this
elliptic curve has self-intersection zero and is a fibre of a holomorphic
fibration of the K3 surface over the Riemann sphere. The
self-intersection number is calculated by the adjunction formula and the
elliptic fibration is defined by a pencil containing the elliptic curve
(see \cite[Propn.~11.1]{hassett}).
All but finitely many fibres of this fibration are smooth and are minimal
submanifolds with respect to a K\"ahler metric on the ambient K3
surface. Moreover, by the Wirtinger theorem the fibres are volume-minimizing
in their homology class.
Harvey and Lawson~\cite{calibrated} generalized the volume-minimizing
property of complex submanifolds of K\"ahler manifolds by developing a
concept of calibrated  minimal submanifolds of Riemannian manifolds. One
type of examples of calibrated submanifolds found in~\cite{calibrated} is
the coassociative submanifolds. These are four-dimensional submanifolds
which occur in the Euclidean $\RE^7$ and, more generally, in 7-dimensional
Riemannian manifolds with holonomy contained in the Lie group $G_2$
(see~\S\ref{synopsis} for precise definitions).

McLean~\cite{mclean} studied the deformation theory for calibrated
submanifolds. In particular, he showed that local deformations of
a smooth compact coassociative submanifold $X$ are always unobstructed
and the corresponding `moduli space' is a smooth manifold of dimension
$b^+(X)$, the positive index of the cup-product on $H^2(X)$.
(This corresponds to the negative index in the actual statement
of~\cite[Theorem~4.5]{mclean} as we use a different sign convention for the
$G_2$-structures.) If $b^+(X)=3$
then, under some reasonable additional assumptions, the local coassociative
deformations of~$X$ give a local foliation. This naturally leads to a
conjecture that some compact 7-manifolds with holonomy $G_2$ are
fibred by coassociative submanifolds, possibly with singular fibres.
Joyce~\cite{joyce-g2,joyce-book} constructed first examples of compact
7-manifolds with holonomy~$G_2$, or $G_2$-manifolds; his results
include examples of compact coassociative calibrated submanifolds
obtained by carefully chosen orientation-reversing involutions. Later
the author gave a series of topologically different examples of compact
$G_2$-manifolds using a `generalized connected sum'
construction~\cite{g2}, see also \cite{k3g2}. It is also shown in~\cite{KN}
that one of the $G_2$-manifolds obtainable as in~\cite{k3g2} is a
deformation of a $G_2$-manifold constructed in~\cite{joyce-g2}. The purpose
of this paper is to construct coassociative fibrations of the
$G_2$-manifolds obtained in~\cite{g2,k3g2}.
A generic fibre of these fibrations is diffeomorphic to
a K3 surface. To the author's knowledge, these are the first examples
of fibrations of {\em compact} smooth manifolds with holonomy $G_2$ by
{\em coassociative minimal} submanifolds.

Further motivation for this paper comes from similarities between
coassociative submanifolds of $G_2$-manifolds and special Lagrangian
submanifolds of Calabi--Yau manifolds. Special Lagrangian submanifolds
are another instance of calibrated submanifolds appearing
in~\cite{calibrated} with unobstructed deformation theory established
in~\cite{mclean}. There are further similarities. For both types of
compact calibrated submanifolds, the `number of moduli' of the local
deformations is given by some Betti number of the submanifold.
Special Lagrangian submanifolds attracted much interest in recent
years in connection with the SYZ conjecture~\cite{SYZ} which explains
mirror symmetry between Calabi--Yau threefolds. The SYZ conjecture was
motivated by studies in string theory and inspired works on the mirror
symmetry for $G_2$-manifolds admitting coassociative fibrations
\cite{ach99,GYZ}. A part of the SYZ conjecture asserts that some
Calabi--Yau threefolds whose complex structures are close to a certain
degenerate limit (the so-called large complex structure limit) admit
fibrations by special Lagrangian tori. This may be compared to the
requirement that in order to construct coassociative fibrations in this
paper we assume that the torsion-free \mbox{$G_2$-structure} on a
compact 7-manifold $M$ corresponds to a point near the
boundary of the moduli space for torsion-free $G_2$-structures.
(This boundary point is represented by a pair of asymptotically cylindrical
manifolds used in the construction of~$M$; this is made precise
in~\cite[\S~5]{nordstrom2}.)
More explicitly, a connected sum $G_2$-manifold should have a
`sufficiently long neck' (as measured by the parameter~$T$ in the Main
Theorem).

On the other hand, there is a difference concerning properties of the
singular fibres. It was shown in~\cite{joyce-syz} that generic fibrations
of Calabi--Yau threefolds by special Lagrangian tori are piece-wise
smooth and only continuous along the loci of the singular fibres which
occur in families of codimension {\em one}.
The coassociative fibrations constructed here are smooth and
are continuously differentiable (more precisely, $C^{1,\alpha}$) at the
singular fibres. The fibres develop singularities modelled on ordinary
double points in the holomorphic 
deformation families of K3 surfaces, and the singular fibres occur in
families of codimension {\em two}, parameterized by a Hopf-type link in
$S^3$. This type of coassociative fibration survives under small
perturbations of the ambient $G_2$-structure, defined by closed 3-forms.

The latter property is crucial as our construction of the coassociative
fibrations relies on perturbative analysis, in the form of implicit
function theorem in Banach spaces. In fact, we build up on the gluing
construction of 
$G_2$-metrics in~\cite{g2} and initially construct an `approximating
fibration' whose fibres are coassociative with respect to a $G_2$-structure
with small torsion. Some familiarity with~\cite{g2} is therefore required
and we include in~\S\ref{pencils} a review of the relevant details before
stating the main result and explaining the strategy of the proof
in~\S\ref{gluing}. The strategy requires, as a preliminary step, an
extension of McLean's deformation theory to coassociative K3 orbifolds
which we carry out in~\S\ref{linear}. Appropriate deformation and stability
result was recently proved by Lotay~\cite{lotay2009}.
With the deformation set-up in place, we proceed to an implicit function
argument for coassociative K3 orbifold fibres in \S\ref{thrm.A} and
complete the construction in~\S\ref{thrm.B} which deals with the smooth
fibres.
\medskip

{\bf Acknowledgements.} 
I am grateful to Simon Donaldson, Dominic Joyce and Miles Reid for helpful
discussions on aspects of this project. I also thank Dominic Joyce for
valuable comments on an earlier version of this paper.

\section{$G_2$-structures and coassociative submanifolds}
\label{synopsis}

We give a short summary of some standard results concerning the
Riemannian geometry associated with the group $G_2$
and coassociative submanifolds associated with the \mbox{$G_2$-structures.}
The results given in this section are mostly gathered
from~\cite{bryant99,calibrated,joyce-book,mclean,salamon} and
the reader is referred there for further details.

The group $G_2$ may be defined as the group of automorphisms of
the cross-product algebra on $\RE^7$ interpreted as the space of pure
imaginary octonions. As the octonions form a normed algebra, any
automorphism in $G_2$ necessarily preserves the Euclidean metric and
orientation of~$\RE^7$. The cross-product can be encoded as an
anti-symmetric 3-linear form $\phi_0$ on $\RE^7$ defined by
$\phi_0(a,b,c)=\langle a\times b,c\rangle$, where
$\langle\cdot,\cdot\rangle$ denotes the Euclidean inner product. The 3-form
$\phi_0$ and its Hodge dual $*\phi_0$ are explicitly expressed as
        \begin{equation}\label{model}
        \begin{split}
\phi_0 &= e_{123} + e_{145} + e_{167} + e_{246} - e_{257} - e_{347} -
e_{356},
\\
*\phi_0 &= e_{4567} + e_{2367} + e_{2345} + e_{1357} - e_{1346} -
e_{1256} - e_{1247},
        \end{split}
        \end{equation}
where $e_{i_1\ldots i_k}=e_{i_1}\wedge\ldots\wedge e_{i_k}$ and
$e_i\in (\RE^7)^*$, $i=1,\ldots,7$ is the standard orthonormal co-frame.
The group $G_2$ may be equivalently defined as the subgroup of
the orientation-preserving linear isomorphisms $GL_+(7,\RE)$ fixing
$\phi_0$ or the subgroup of $GL_+(7,\RE)$ fixing $*\phi_0$, in the
action on, respectively, 3- or 4-forms. The group $G_2$ is a compact
Lie group of dimension~14 and thus a proper subgroup of~$SO(7)$.

Let $M$ be a smooth oriented 7-manifold. Denote by $\Omega^3_+(M)$ the
subset of 3-forms $\phi$ on~$M$ whose value $\phi(p)$ at any point
$p\in M$ can be identified with $\phi_0$ via some orientation-preserving
linear isomorphism $T_pM\to\RE^7$. Any $\phi\in\Omega^3_+(M)$ induces
a $G_2$-structure on~$M$; we shall sometimes, slightly inaccurately, say
that $\phi$ {\em is} a $G_2$-structure. The $GL(7,\RE)$-orbit of $\phi_0$
is open in $\Lambda^3(\RE^7)^*$, so the subset
$\Omega^3_+(M)\subset\Omega^3(M)$ is open in the uniform norm 
topology. As $G_2\subset SO(7)$ any $\phi\in\Omega^3_+(M)$ induces a
metric $g(\phi)$ and an orientation of~$M$, so that the image of
$e_1,\ldots,e_7$ via the isomorphisms $T_pM\to\RE^7$ identifying $\phi_p$
with $\phi_0$ is an orthonormal positive co-frame at each $p\in M$.  The
4-form $*_\phi\phi$ on $M$ defined using the Hodge star of the metric
$g(\phi)$ is, in a similar sense, point-wise 
isomorphic to $*\phi_0$ in~\eqref{model}.

The metric $g(\phi)$ induced by a $G_2$ structure $\phi\in\Omega^3_+(M)$
will have holonomy contained in~$G_2$ if and only if
        \begin{equation}\label{g2torsionfree}
d\phi=0,\qquad d*_\phi\!\phi=0,
        \end{equation}
the $G_2$-structure $\phi$ then is said to be {\em torsion-free}. 
Conversely, if the holonomy of a metric $g$ is contained in $G_2$ then
there is a $\phi\in\Omega^3_+(M)$ satisfying~\eqref{g2torsionfree}
so that $g=g(\phi)$. \cite[Lemma 11.5]{salamon}  The pair $(M,\phi)$ then
is called a {\em $G_2$-manifold} and a metric $g(\phi)$ induced by a
torsion-free $G_2$-structure on a 7-manifold $M$ is called a
{\em $G_2$-metric}. On a compact $M$, the holonomy of a $G_2$-metric is
exactly~$G_2$ if and only if the fundamental group of~$M$ is
finite~\cite[Propn.~10.2.2]{joyce-book}.
\medskip

Now we recall the concept of coassociative submanifolds.
Any 4-tuple of orthogonal unit vectors $v_1,v_2,v_3,v_4\in\RE^7$ satisfies
the inequality
        \begin{equation}\label{ca.ineq}
*\phi_0(v_1,v_2,v_3,v_4)\le 1.
        \end{equation}
If the equality in attained in~\eqref{ca.ineq} then the oriented
4-dimensional subspace of~$\RE^7$ defined by $v_1,v_2,v_3,v_4$ is called a
{\em coassociative subspace}. More generally, if $M$ is a 7-manifold
and $\phi\in\Omega^3_+(M)$ is a $G_2$-structure form then any oriented
4-dimensional submanifold $X\subset M$ satisfies 
        \begin{equation}\label{ca.gen}
*_\phi\phi|_X\le \vol_X,
        \end{equation}
in the sense that $*_\phi\phi|_X=\alpha \vol_X$ for some real constant
$\alpha\le 1$, where $\vol_X$ is the volume form of the metric on~$X$
induced by the embedding.
        \begin{prop}[cf.~\mbox{\cite[Cor.~IV.1.20]{calibrated}}]
For an orientable 4-dimensional submanifold~$X$ of a \mbox{7-manifold} $M$
with a $G_2$-structure $\phi\in\Omega^3_+(M)$, the equality
$*_\phi\phi|_X=\vol_X$ is attained in~\eqref{ca.gen} for some (necessarily
unique) orientation of~$X$ if and only if $\phi|_X=0$.
        \end{prop}
If a $G_2$-structure on~$M$ satisfies $d\!*_\phi\!\phi=0$
then the 4-form $*_\phi\phi$ is an instance of a calibration, called
{\em coassociative calibration}, and a 4-dimensional submanifold $X$ with
$*_\phi\phi|_X=\vol_X$ is then said to be calibrated by $*_\phi\phi$.
Any compact calibrated submanifold minimizes the volume in its
homology class, in particular a compact $X$ calibrated by $*_\phi\phi$
is a minimal submanifold of the Riemannian manifold~$(M,g(\phi))$.

The $G_2$ structures that we deal with in this paper will always be
defined by {\em closed} \mbox{3-forms}. The following terminology
will be in use. A 4-dimensional submanifold $X\subset M$ will be called
{\em coassociative}, with respect to a $G_2$-structure
$\phi\in\Omega^3_+(M)$, if $\phi|_X=0$.  If, in
addition, this $G_2$-structure satisfies $d*_\phi\phi=0$ (and so
$*_\phi\phi$ is a calibration) then we say that $X$ is a
{\em coassociative calibrated} submanifold.

A remarkable property of compact coassociative submanifolds is that
their local deformations are unobstructed, a result due to
McLean~\cite{mclean}. It uses a `tubular neighbourhood theorem'
(e.g.~\cite[Ch.~IV]{lang}), which may be stated as follows.
	\begin{thrm}\label{tubular}
For any embedded submanifold~$X$ of a Riemannian manifold $(M,g_M)$, the
exponential map of $g_M$,
$$
v(p)\in (N_{X\subset M})_p\subset T_pM\mapsto \exp_{v(p)}\in M,
\qquad
p\in X,
$$
induces a diffeomorphism of a neighbourhood
$$
\bigl\{v(p)\in N_{X/M}:|v|<C_X\min\{\delta(p),|II(p)|\}
\text{ if }v(p)\in T_pM,\break  p\in X\bigr\}
$$
of the zero section of the normal bundle $N_{X/M}$ onto a neighbourhood 
of $X$ in~$M$. Here $\delta(p)$ denotes the injectivity radius of $g_M$
and $II(p)$ the second fundamental form at~$p$, $C_X>0$ is a constant
depending only on~$X$.
	\end{thrm}
As $\delta(p)$ depends continuously on~$p\in X$ it attains a maximum
$\delta>0$ when $X$ is compact. In that case, the inequality in
Theorem~\ref{tubular} defining a tubular neighbourhood can be simplified
to $|v(p)|<\delta$.

Thus any local deformation of compact coassociative~$X$ can be
written as $X_v=\exp_v(X)$, for some $v\in\Gamma(N_{X\subset M})$ with
small $\|v\|_{C^0}$ in the metric induced from~$M$. Suppose that
$\phi|_X=0$. A local deformation $X_v$ of~$X$ will be coassociative if and
only if $F(v)=0$, where
        \begin{equation}\label{mclean}
F: v\in\Gamma(N_{X\subset M})\to F(v)=(\exp_v)^*\phi\in\Omega^3(X),
        \end{equation}
The normal bundle $N_{X\subset M}$ of a coassociative submanifold is
canonically isomorphic and isometric to the bundle of self-dual 2-forms
$\Lambda^+T^*X$ via
        \begin{equation}\label{normal}
v\in\Gamma(N_{X/M})\to v\lrcorner\phi|_X \in\Omega^+(X). 
        \end{equation}
(The $G_2$ structure 3-forms $\phi$ used in~\cite{mclean} differ
from~\eqref{model} by the opposite sign, which leads to $\Omega^-(X)$ rather
than~$\Omega^+(X)$ in the right-hand side of~\eqref{normal}.)
Composing~$F$ with the inverse of~\eqref{normal}, we obtain a map
$$
\check{F}:\Omega^+(X)\to\Omega^3(X).
$$
McLean proves:
          \begin{thrm}[cf.~\cite{mclean}, \S 4]\label{CA}
Suppose that a $G_2$-structure on a 7-manifold~$M$ is given by a closed
3-form and $X\subset M$ is a coassociative submanifold. Then:

(a) the derivative $(d\check{F})_0$ is given by the exterior
derivative $d:\Omega^+(X)\to\Omega^3(X)$, and

(b) the image of $\check{F}$ consists of exact forms.
          \end{thrm}
          \begin{remarks}
The statement of Theorem~\ref{CA} extracts a part of McLean's results
which does not require the compactness of~$X$.
The result is stated in~\cite{mclean} for coassociative calibrated
submanifolds but it was later observed in~\cite{goldstein1} that the
condition $d*_\phi\phi=0$ is not used in the proof.
          \end{remarks}
If $X$ is compact then any exact 3-form on a compact oriented
Riemannian 4-manifold is the differential of a self-dual form and
Theorem~\ref{CA} sets the scene for an application of the implicit function
theorem in Banach spaces. It follows that the local deformations of a
compact coassociative $X$ form a smooth manifold of dimension equal to the
dimension $b^+(X)$ of harmonic self-dual forms on~$X$, see
\cite[Theorem 4.5]{mclean} or \cite[Theorem 2.5]{joyce-salur}.

Bryant~\cite{bryant99} proved that any closed real-analytic oriented
Riemannian 4-manifold $X$ with trivial $\Lambda^+T^*X$ arises as a
coassociative calibrated submanifold in some manifold with torsion-free
$G_2$-structure. If there are three harmonic self-dual
forms on~$X$ that are linearly independent at every point then
Bryant's result produces examples of non-compact (local) $G_2$-manifolds
foliated by the coassociative deformations of~$X$.

Finally, an application of McLean's theory shows that compact
coassociative submanifolds are `stable' under small deformations of
the $G_2$ structure. The property will be crucial for the main results
of this paper, the following theorem provides an introduction.
        \begin{thrm}[cf.~\mbox{\cite[Theorem 12.3.6]{joyce-2007}}]\label{st}
Suppose that $\phi(s)\in\Omega^3(M)$, $s\in\RE$, is a smooth path of
closed $G_2$-structure forms on~$M$, and $X$ is a compact submanifold of~$M$
such that $\phi(0)|_X=0$ and the form $\phi(s)|_X$ is exact for any~$s$.
Then there there is an $\eps>0$ and for each $|s|<\eps$ a section
$v(s)$ of $N_{X/M}$ smoothly depending on~$s$, such that $v(0)=0$
and $\phi(s)$ vanishes on $\exp_{v(s)}(X)$.
        \end{thrm}

It is worth to point out the different roles of the two equations
in~\eqref{g2torsionfree}: 
the second equation ensures that compact coassociative submanifolds
are minimal, whereas the first equation ensures that they have a good
deformation theory (Theorems~\ref{CA} and~\ref{st}).
Moreover, it is noted in~\cite[0.3.3]{bryant99} that the generic
$G_2$-structure satisfying $d*_\phi\phi=0$ (but not necessarily
$d\phi=0$) will not admit any coassociative submanifolds.

\section{From pencils of K3 surfaces to the approximating fibrations}
\label{pencils}

The method of construction of compact irreducible $G_2$ manifolds that we
shall consider was developed in~\cite{g2}. Because some technical details
of this construction will be important in what follows we shall review these
details here. We shall also deduce some immediate consequences
concerning coassociative submanifolds which are not given in~\cite{g2}.

Let $W$ be a Ricci-flat K\"ahler complex threefold and $\omega$ the
K\"ahler form on $W$. If the holonomy of the K\"ahler metric is
contained in $SU(3)$ (which will be the case e.g.\ if $W$ is
simply-connected) then there is a nowhere-vanishing holomorphic
$(3,0)$-form $\Omega$ on~$W$, sometimes called a holomorphic volume form.
We shall call the pair $(\omega,\Omega)$ a {\em Calabi--Yau structure}
on~$W$. A Calabi--Yau structure $(\omega,\Omega)$ induces a torsion-free
$G_2$-structure on the 7-manifold $W\times S^1$ given by the 3-form
        \begin{equation}\label{product}
\phiCY = \omega\wedge d\theta + \im\Omega,
        \end{equation}
where $\theta$ is the standard `angle coordinate' on $S^1$
(cf.~\cite[Propn.~11.1.2]{joyce-book}, our holomorphic volume form differs
from the one used there by the factor~$i$). The form $\phiCY$ induces
the product metric $g_W+d\theta^2$ on $W\times S^1$ corresponding to the
Ricci-flat K\"ahler metric $g_W$ on~$W$ and we shall sometimes refer to
$\phiCY$ as a product $G_2$-structure.

If $X$ is a complex surface in~$W$ then, for any $\theta_0\in S^1$,
$$
\phiCY|_{X\times\{\theta_0\}}=\im(\Omega|_X)=0,
$$
so $X\times\{\theta_0\}$ is a coassociative calibrated
submanifold of $W\times S^1$. \cite{mclean,joyce-book}

Let $V$ be a (non-singular) compact complex threefold with $c_1(V)>0$,
i.e.\ a Fano threefold,
and $D\in|-K_V|$ a K3 surface in the anticanonical linear system
of~$V$. Let $D'\in|-K_V|$, $D'\neq D$, be another K3 surface in the
anticanonical class, so that $C=D'\cap D$ is a non-singular connected curve
in~$V$. Blowing up~$C$ we obtain a new threefold $\tV$ with a holomorphic
map $\tau:\tV\to\CP^1$ whose fibres are proper transforms of the surfaces
in the pencil defined by $D$ and~$D'$.

The proper transform $\tD$ of~$D$ is an anticanonical divisor on $\tV$
and the complement non-compact complex threefold $W=\tV\minus\tD$ has
trivial canonical bundle. We can define
a holomorphic coordinate, $\zeta$ say, on $\CP^1$, so that $D=\tau^{-1}(0)$.
The fibre $\tau^{-1}(\zeta)$ is diffeomorphic to $\tD$ as a real
4-manifold whenever $|\zeta|$ is sufficiently small, thus the real
6-manifold underlying $W$ has a cylindrical end
$\tau^{-1}(\{0<|\zeta|<\eps\})$ diffeomorphic to
$\RE_{>0}\times S^1\times D$. Denote the real coordinates on the first two
factors by $t,\theta$ and then $\zeta=e^{-t-i\theta}$.
The complex structure on the end of $W$ is asymptotic, as
$t\to\infty$, to the product complex structure
on~$\RE_{>0}\times S^1\times D$, where $D$ is considered
with the complex structure and K\"ahler form induced by the embedding in~$V$.
Note that the latter condition on the K\"ahler form is {\em not}
restrictive: every K\"ahler metric on~$D$ is obtainable as a restriction of
some K\"ahler metric on~$V$. \cite{g2} \cite{k3g2}

Another class of complex threefolds with similar properties was recently
constructed in~\cite{k3g2}. It uses K3 surfaces with non-symplectic
involution.

Recall that by Yau's solution of the Calabi conjecture~\cite{yau} the K\"ahler
K3 surface $D\subset V$ has a unique Ricci-flat K\"ahler metric in its
K\"ahler class. We shall write $\kappa_I$ for the K\"ahler form of
this metric and $\kappa_J+i\kappa_K$ for a holomorphic volume form on~$D$.
The following `non-compact version of the Calabi conjecture' for~$W$
is proved in~\cite{g2} using the results of~\cite{TY1}.
        \begin{thrm}[\cite{g2},  \S\S$\,$2--3,$\,$6 ]\label{bCY}
The threefold~$W$ is simply-connected and has a complete Ricci-flat
K\"ahler metric $g_W$ of holonomy $SU(3)$. This metric is
asymptotically cylindrical in the sense that on the cylindrical end
$\tau^{-1}\{0<|\zeta|<\eps\}\subset W$ the K\"ahler form $\omega$ of $g_W$ has
an asymptotic expression
$$
\omega|_{\tau^{-1}\{0<|\zeta|<\eps\}}=dt\we d\theta+\kappa_I+d\psi,
$$
and there is a holomorphic volume form $\Omega$ on~$W$ with an asymptotic
expression
$$
\Omega|_{\tau^{-1}\{0<|\zeta|<\eps\}}=
(dt+id\theta)\we(\kappa_J+i\kappa_K)+d\Psi,
$$
where $\zeta=e^{-t-i\theta}$ and the differential forms
$\psi,\Psi$ and all their derivatives decay at the rate $O(e^{-\lambda_W t})$
along the end of~$W$, with the exponent $\lambda_W>0$ depending only on the
K\"ahler metric $\kappa_I$ on~$D$.
        \end{thrm}
For any $W$ satisfying the assertion of Theorem~\ref{bCY}, the
7-manifold $W\times S^1$ has a cylindrical end
$\RE_{>0}\times S^1\times D\times S^1$ and the product $G_2$-structure
$\phiCY$ on $W\times S^1$ has an asymptotic expression
$$
\phiCY|_{\tau^{-1}\{0<|\zeta|<\eps\}\times S^1}=\phi_D+ d\tilde\psi,
$$
where
\begin{align}
&\phi_D = dt\we d\theta\we d\theta' + \kappa_I\we d\theta' +
\kappa_J\we d\theta + \kappa_K\we dt,\label{su2g2}\\
&\tilde\psi=\psi\we d\theta'+\im\Psi,\notag
\end{align}
and $d\theta'$ denotes the standard non-vanishing 1-form on the last
$S^1$ factor.
Note that the asymptotic model $\phi_D$ is determined by the
choice of the K3 divisor $D\in |-K_V|$ alone and does not depend on the
choice of $D'$ and the resulting pencil on~$V$.
        \begin{prop}\label{genpen}
For a generic choice of Fano threefold~$V$ in its deformation family
and a generic choice of $D'\in|-K_V|$, the fibres of the map
$\tau:\tV\to\CP^1$ define a `generic  Lefschetz fibration' in the following
sense:\\
(1) the critical points of $\tau$ are non-degenerate (Morse points): if
$d\tau(w)=0$ then the Hessian of $\tau$ at $w$ is
non-singular, and\\
(2) any fibre $\tau^{-1}(\zeta)$ contains at most one critical point of~$\tau$.
        \end{prop}
        \begin{pf}
This is an application of~\cite{g2} and some standard results in algebraic
geometry and we only give an outline of the proof. A generic anticanonical
divisor on~$V$ is a smooth surface of type~K3~\cite{sh1}. The K3 surfaces
arising as smooth anticanonical divisors on the deformations of $V$ form
a Zariski open subset in a moduli space $\eus{K}$ of K3 surfaces whose Picard
lattice contains a copy of $H^2(V,\ZE)$ as a sublattice~\cite[\S 7]{g2}.
The moduli space $\eus{K}$ is a quasiprojective complex algebraic variety
of dimension $20-b^2(V)$. The degenerations of the K3 surfaces in $\eus{K}$
developing an ordinary double point are generic.
The anticanonical linear system $|-K_V|$ is parameterized by $\CP^N$
and an application of the Riemann--Roch theorem gives the dimension
$N=-K_V^3/2+2$, so $N\ge 3$ \cite[Propn.~1.3]{Is}. 

The singular anticanonical divisors on~$V$ are therefore parameterized
by an algebraic subvariety $S$ of $\CP^N$ of codimension at least 1
and one can show that for a generic Fano threefold in the deformation
family of~$V$ any connected component of $S$ contains a K3 orbifold with
the only singularity an ordinary double point. Each
of the conditions (1) and (2) in Proposition~\ref{genpen} is an open 
condition in the Zariski topology of~$S$. Violation
of (1) or (2) defines a further subvariety of positive
codimension in~$S$. This latter subvariety therefore has codimension
at least~2 in $\CP^N$ and can be avoided in the pencil
through $D$ and~$D'$ with a generic choice of~$D'$.
        \end{pf}

If $w_0\in W$ is a critical point of~$\tau$ and the Hessian of $\tau$
at~$w_0$ is non-degenerate then by the Morse lemma there is a system of local
holomorphic coordinates $z_j$ near~$w_0$, such that
$\tau=z_1^2+z_2^2+z_3^2$ in these coordinates.
The fibre $X_0$ through $w_0$ is an orbifold with an isolated singularity
$z_1^2+z_2^2+z_3^2=0$ at~$w_0$. (For a general theory of orbifolds
see~\cite{baily}.) A neighbourhood of~$w_0$ in $X_0$ admits a local
parameterization by uniformizing coordinates $u_1,u_2\in\CX^2/\pm1$,
        \begin{equation}\label{unicoord}
(u_1,u_2)\mapsto\bigl(i(u_1^2+u_2^2),\; u_1^2-u_2^2,\; 2u_1u_2\bigr),
        \end{equation}
so this neighbourhood is homeomorphic to a cone on~$\RP^3$.

A K3 fibration $\tau$ has only finitely many singular fibres. Every
singular fibre $X_0$ of a generic $\tau$ (in the sense of
Proposition~\ref{genpen}) is an `orbifold K3 surface' with unique
singularity which is an ordinary double point and $\tau$ near $X_0$
defines a one-parameter family of non-singular deformations of~$X_0$.
On the other hand, it is well-known (e.g.~\cite[\S 2.6]{aspinwall}) that by
{\em blowing up} $w_0$ one achieves a resolution of the $\pm1$ orbifold
singularity of~$X_0$ which is again a K3 surface (in general,
not isomorphic to any nearby fibre of~$\tau$) and the exceptional divisor
is a complex curve with self-intersection~$-2$.

The number $\mu_V$ of singular fibres of a generic $\tau$ is the number of
K3 orbifolds in the corresponding generic Lefschetz pencil on~$V$. This
number is calculated by an application of Lefschetz theory of hyperplane
sections to the topology of algebraic varieties~\cite[\S 5]{AF},
        \begin{equation}\label{mu}
\mu_V=2\chi(D)-\chi(C)-\chi(V)=48+(-K_V^3)-\chi(V)
        \end{equation}
where $C$ is the blow-up locus (the axis of the pencil of K3
surfaces in~$V$), so $-\chi(C)=-K_V^3$ for a Fano threefold~$V$
\cite[Propn.~1.6]{Is}.
\smallskip

We may assume, rescaling $\zeta$ by a non-zero constant factor if
necessary, that there are no singular fibres of~$\tau$ on the cylindrical 
end $\{t>0\}\subset W$, where $t=-\log|z|$ as before. Extend~$t$ to
a smooth function, still denoted by~$t$, defined on all of~$W$ with
$t<0$ away from the cylindrical end. Fix once and for all a smooth cut-off
function $\alpha(s)$ with $\alpha(s)=0$, for $s\le 0$, and $\alpha(s)=1$,
for $s\ge 1$. The 3-form
        \begin{equation}\label{phiT}
\phi_{W,T}=\phiCY-d(\alpha(t-T)\tilde\psi)
        \end{equation}
defines a $G_2$-structure on $W_T\times S^1$ for every sufficiently
large~$T$. It interpolates between the product torsion-free
$G_2$-structure $\phiCY$ on $W\times S^1$ induced by the
Ricci-flat K\"ahler structure on~$W$ and the product torsion-free
$G_2$-structure $\phi_D$ on the half-cylinder \mbox{$[T-1,\infty)\times
S^1\times D\times S^1$} induced by the Calabi--Yau (hyper-K\"ahler)
structure forms $\kappa_I$, $\kappa_J,\kappa_K$ on~$D$ as in~\eqref{su2g2}.
        \begin{prop}\label{in}
For any $T>0$, the 3-form $\phi_T$ vanishes on each fibre of the map
$\tau\times\id_{S^1}:W\times S^1\to\CP^1\times S^1$.
        \end{prop}
        \begin{pf}
It is clear that the claim is true away from the cut-off region
$R=[T-1,T]\times S^1\times D\times S^1$ because the map $\tau$ is
holomorphic both in the complex structure of~$W$ and in the product
complex structure of $\RE\times S^1\times D$ on the end of~$W$.

On the cut-off region, we have
$$
\phi_T|_{R}=(1-\alpha)\phiCY+\alpha\phiD-\alpha'dt\we\tilde\psi
$$
and the claim follows as any fibre of $\tau$ on the end of $W\times S^1$
is contained in a level set $\{t=\const\}$.
        \end{pf}

Now let $W_1$ and $W_2$ be two asymptotically cylindrical Calabi--Yau
threefolds given by Theorem~\ref{bCY}
and define $W_j(T)=W_j\minus\{t_j>T\}$, for any $T>2$, $j=1,2$.
Assume that the respective two hyper-K\"ahler K3 surfaces $D_j$ are
`hyper-K\"ahler rotations' of each other which means that there is an
isometry $f:D_1\to D_2$ of the Riemannian 4-manifolds such that $f^*$
interchanges the K\"ahler and holomorphic volume forms of $D_1$ and $D_2$ as
follows: $\kappa^{(2)}_I\mapsto\kappa^{(1)}_J$,
$\kappa^{(2)}_J\mapsto\kappa^{(1)}_I$, 
$\kappa^{(2)}_K\mapsto-\kappa^{(1)}_K$.
(Such an isometry $f$ always exists after some deformations of the Fano
threefolds $V_1,V_2$ used in the construction of
$W_1,W_2$~\cite[Theorem~6.44]{g2}. See also \cite[Theorem~5.3]{k3g2})

Construct a compact 7-manifold $M$ by joining $W_1(T)$ and
$W_2(T)$
        \begin{subequations}\label{gsum}
        \begin{equation}
M=\bigl(W_1(T)\times S^1\bigr)\cup_\Upsilon \bigl(W_2(T)\times S^1\bigr),
        \end{equation}
identifying collar neighbourhoods of the boundaries via the
orientation-preserving diffeomorphism
        \begin{multline}
\Upsilon:
(y,\theta_1,\theta_2,T+t)\in D_1\times S^1\times S^1\times]T+1,T+2[\to\\
(f(y),\theta_2,\theta_1,T+3-t)\in D_2\times S^1\times S^1\times ]T+1,T+2[
        \end{multline}
        \end{subequations}
Then $\Upsilon^*\phi_{D_2}=\phi_{D_1}$, so the two $G_2$-structures
$\phi_{j,T}\in\Omega^3_+(W_j(T)\times S^1)$ defined by~\eqref{phiT}
agree on the overlap and together give a
well-defined $G_2$-structure $\phi_T\in\Omega^3_+(M)$
on the compact 7-manifold. This latter $G_2$-structure form satisfies
	\begin{equation}\label{errterm}
d\phi_T=0,\qquad \|d*_{\phi_T}\phi_T\|_{L^p_k}
<C_{p,k}e^{-\lambda T},
\end{equation}
for each $T>2$, where $\lambda<\lambda_{W_j}$, $j=1,2$, and $\lambda_{W_j}$
are defined in Theorem~\ref{bCY}.
        \begin{thrm}[\cite{g2}, \S5]\label{glug2}
There exists $T_0$ and for each $T>T_0$ a 2-form $\eta=\eta_T$ on the
7-manifold~$M$ satisfying $\|\eta_T\|_{L^p_k}<K_{p,k}e^{-\lambda T}$
and such that
$$
d*(\phi_T+d\eta_T)=0,
$$
with the Hodge star $*$ in the above formula defined by the metric
$g(\phi_T+d\eta_T)$. Thus the 3-form $\phi_T+d\eta_T$
induces a torsion-free $G_2$-structure on~$M$ and the holonomy group of
the metric $g(\phi_T+d\eta_T)$ is~$G_2$.
        \end{thrm}
The last claim is true because $M$ is simply-connected.
The 3-form $\phi_T$ may be thought of as an approximation, improving
as $T\to\infty$, of a torsion-free $G_2$-structure on~$M$.

Return, for the moment, to the $G_2$-structures $\phi_{j,T}$
on the two pieces of~$M$. It is clear from the construction~\eqref{gsum}
of~$M$ that the fibrations $\tau^{(j)}$ of $W_j(T)\times S^1$ agree on the
overlap and hence can be patched to define, for any $T>1$ a fibration
$\tau_T$ of~$M$ as shown in the following commutative diagram
        \begin{equation}\label{solid.torus}
        \begin{CD}
W_{1,T}@. \times S^1 \;\cup\; @. W_{2,T}@. \times S^1 @= M\\
@V\tau^{(1)}VV @. @V\tau^{(2)}VV @. @VV\tau_T V \\
\Delta@. \times S^1 \;\cup\; @. \Delta@. \times S^1 @= S^3 .
        \end{CD}
        \end{equation}
where the bottom row is a well-known splitting of the 3-sphere
into two solid tori ($\Delta$ denotes a disc in $\CP^1$). We obtain
from~\eqref{solid.torus} and Proposition~\ref{in} the main result of this
section.
        \begin{thrm}\label{approx}
Let $(M,\phi_T)$ be a compact 7-manifold with a $G_2$-structure constructed
from a pair of Fano threefolds, as defined above.
Then $\tau_T:M\to S^3$ is a coassociative fibration defined
by~\eqref{solid.torus}, with respect to $\phi_T$. The singular fibres
of~$\tau_T$ form a subset of codimension~2 in~$M$ and are projected
by~$\tau$ onto a link in~$S^3$.

This link consists of $\mu_{V_1}$ (disjoint) circles in
$S^3\minus\tau_T(W_1(T)\times S^1)$ and $\mu_{V_2}$ (disjoint) circles in
$S^3\minus\tau_T(W_2(T)\times S^1)$, where $\mu_{V_j}$ is defined
in~\eqref{mu}. Each circle in $S^3\minus\tau_T(W_j(T)\times S^1)$
is linked, with linking number $1$, with each circle in the
other subset $S^3\minus\tau_T(W_{3-j}(T)\times S^1)$
and is not linked with any circle in $S^3\minus\tau_T(W_j(T)\times S^1)$.
        \end{thrm}
Of course, the 4-form $*_{\phi_T}\phi_T$ is not in general a calibration
on~$M$ and the fibres of $\tau_T$ are not necessarily calibrated by
$*_{\phi_T+d\eta_T}(\phi_T+d\eta_T)$. 
The estimate on the 2-form $\eta_T$ given in Theorem~\ref{glug2} yields an
upper bound on $d\eta_T|_{X}$, for each non-singular fibre $X$ of $\tau_T$,
        \begin{equation}
\|d\eta_T|_{X}\|_{L^p_k} < \tilde{K}_{p,k}(X) e^{-\lambda T},
        \end{equation}
for each $T>T_0$. The constant $\tilde{K}_{p,k}(X)$ in~\eqref{errterm}
depends on a particular choice of norm (more precisely, on the value of
$k-4/p$) as well as on the choice of $X$. We shall return to this
later, in Proposition~\ref{et}.

\section{The gluing theorem for coassociative K3 fibrations}
\label{gluing}

The estimate~\eqref{errterm} on the failure of the fibres of $\tau_T$ to be
calibrated by the 4-form $*_{\phi_T+d\eta_T}(\phi_T+d\eta_T)$ also suggests
that, the map $\tau_T$ might be in some sense an approximation of a
fibration of the holonomy-$G_2$ manifold $(M,g(\phi_T+d\eta_T))$ with
coassociative calibrated fibres. The next theorem---which is the main
result of this paper---asserts that this is indeed the case.
        \begin{varthrm}[Main Theorem]
Let $(M,g(\phi_T+d\eta_T))$ be a compact 7-manifold with holonomy~$G_2$
constructed from a pair of Fano threefolds as defined in the previous
section and let $\tau_T$ be the fibration of~$M$ given by
Theorem~\ref{approx}. There exists, for every sufficiently large~$T$, a
diffeomorphism $h_T$ of $M$ onto itself, exponentially close to $\id_M$,
$$
\sup_{x\in M}\;\dist_{g(\phi_T)}(h_T(x),x)<\const\cdot e^{-\lambda T}
$$
with $\lambda>0$ as determined in Theorem~\ref{glug2},
and such that the fibres of $\tau_T\circ h_T:M\to S^3$ are
coassociative calibrated by $*_{\phi_T+d\eta_T}(\phi_T+d\eta_T)$.
In particular, smooth fibres of $\tau_T\circ h_T$ are minimal submanifolds
of $(M,g(\phi_T+d\eta_T))$.
The map $h_T^{-1}$ can be taken to be $C^1$ on the locus of the singular
fibres of $\tau$ and $C^\infty$ elsewhere on~$M$.

A generic fibre of $\tau_T\circ h_T$ is diffeomorphic to the real
4-manifold underlying a K3 surface. Each singular fibre is an orbifold
diffeomorphic to a `K3 orbifold surface' with one ordinary double point
and no other singularities. The discriminant locus (image of the singular 
fibres) of $\tau_T\circ h_T$ is a link in~$S^3$, as described in
Theorem~\ref{approx}.
        \end{varthrm}
        \begin{remarks}
(1) The smooth K3 fibres of~$\tau_T\circ h_T$ have $b^+=3$ and thus form a
maximal deformation family, by McLean's results~\cite{mclean}. We shall
see, after some additional work below, that the singular fibres
of~$\tau_T\circ h_T$ also form a maximal deformation family.

Another well-known compact 4-manifold with $b^+=3$ is a 4-torus. A
fibration by coassociative 4-tori was obtained by 
Goldstein~\cite{goldstein1} for a compact 7-manifold constructed
in~\cite{joyce-g2}. However, the $G_2$-structure used 
in~\cite{goldstein1} only has a closed 3-form but does {\em not} define
a coassociative calibration (it is close to a calibration and one can
obtain an estimate similar to~\eqref{errterm}), so the fibres need
not be minimal submanifolds.  The singular coassociative fibres
in~\cite{goldstein1} have non-isolated singularities and it appears that
the problem of perturbing into a map with calibrated fibres would require a
different analytic technique than that developed for the fibration~\eqref{solid.torus}. 
We hope to return to this problem in a future paper.

(2) A K3 fibre $X$ of~$\tau$ has trivial normal bundle in~$M$ and hence 
trivial bundle of self-dual forms $\Lambda^+T^*X$. A trivialization of
$\Lambda^+T^*X$ induces an $Sp(1)$-structure on $X$ and hence
a triple of orthogonal almost complex structures $I,J,K$ (relative to
the metric on~$X$) satisfying the quaternionic relations $IJ=-JI=K$. From
the construction of the holonomy-$G_2$ metric on~$M$ one can see that one
of $I,J,K$ is a deformation of the complex structure induced by embedding
of~$X$ in the Calabi--Yau threefold $W_j$. As $T$ tends to infinity, the
metric induced by $g(\phi_T)$ on each fibre of $\tau_T$ converges uniformly
with all derivatives to a K\"ahler metric. For fibres near the middle of
the neck of~$M$, the limit metric is, moreover, hyper-K\"ahler. However,
there is no general reason for a holonomy reduction for the induced metric
on the fibres, for any finite~$T$.
        \end{remarks}
In the rest of the paper we prove the Main Theorem.

We shall construct $h_T$ in the form $h_T=\exp_{v_T}$, for a
smooth vector field $v_T$ on~$M$ (more precisely, $h_T$ will be obtained as
a composition of two exponential maps). A map $\exp_v$ of the compact Riemannian
manifold $(M,g(\phi_T))$ is well-defined whenever the uniform norm of~$v$ is
less than the injectivity radius of~$M$. Furthermore, $\exp_v$ defines a
{\em diffeomorphism} of~$M$ whenever the uniform norm of both $v$ and its first
derivatives is sufficiently small, so that $\exp_v$ is a local
diffeomorphism near every point and a bijection of~$M$.
        \begin{prop}
Let $M$ be a compact 7-manifold with a one-parameter family of Riemannian
metrics induced by the $G_2$-structures $\phi_T\in\Omega^3_+(M)$ defined
by the connected sum construction, as in~\S\ref{pencils}.

Then there exists $\eps>0$ so that for any metric $g(\phi_T)$, $T\ge 1$,
the map $\exp_{v}$ is a diffeomorphism of~$M$ onto itself whenever 
$\|v\|_{C^1}<\eps$.
        \end{prop}
        \begin{pf}
The upper bound $\eps$ can be determined by considering restrictions of
the metric $g(\phi_T)$ to small neighbourhoods and taking the supremum.
The asymptotically cylindrical properties of the metrics constructed on
$W_j\times S^1$ imply that the $\eps$ can be taken positive on these
manifolds and hence also on~$M$ as $g(\phi_T)$ is exponentially close to
the asymptotically cylindrical metrics as $T\to\infty$.
        \end{pf}

We require a $C^1$-small vector field~$v=v_T$ on~$M$ satisfying
$$
\exp_v^*(\phi_T+d\eta_T)|_X=0
$$
for each fibre~$X$ of $\tau_T$, given that $\phi_T|_X=0$ and $\eta_T$ can
be taken as small as we like by choosing a large~$T$. The construction of
such $v$ can be thought of as an infinite-dimensional version of the
implicit function problem $F_X(v,s)=0$ on each fibre $X$, for a family of
vector fields $v=v(s)$, $0\le s\le 1$, with $v(0)=0$, where
$$
F:(v,s)\in\Gamma(TM|_{X})\times[0,1]\to 
\exp_v^*(\phi+s\,d\eta)|_{X}\in\Omega^3(X)
$$
is a variant of the map appearing in McLean's theory (cf.~\S\ref{synopsis}).
Here we temporarily dropped the dependence on~$T$ from the notation.

If all the coassociative fibres of~$\tau$ were smooth, then the desired
vector field $v$ would be easily obtained by a slight modification of
Theorem~\ref{st} for local deformation families of coassociative
submanifolds and then patching together finitely many of these local families,
using the compactness of~$M$. However, the fibres of~$\tau$ develop
singularities. Recall also that the deformation problem for a coassociative
submanifold is expressed as an equation for sections of vector bundles on
the actual submanifold. In light of this, the proof of the Main Theorem
naturally falls into two parts concerned, respectively, with the singular
fibres of~$\tau_T$ and the nearby smooth fibres with `large' curvature.

In order to implement the implicit function strategy
we require, in the first place, an
extension of the deformation theory to the coassociative {\em K3 orbifolds}
arising as the singular fibres. More precisely, the required property concerns
the linearization $(D_1F)_0$ of the deformation map in~$v$ at $v=0$. This map
should be a surjective Fredholm map between appropriate Banach spaces
(weighted Sobolev spaces in Theorem~\ref{noobstr}), so the deformations of
coassociative K3 orbifold fibres of $\tau_T$ are {\em unobstructed}.

The perturbation~$h_T$ of $\tau_T$ will be obtained as a
composition of exponential maps via the following two results.
The proof of Theorem~A requires a `stability' result for coassociative
cones defined by complex 2-dimensional cones in~$\CX^3$. Appropriate
stability result for the tangent cones at the singular points
of K3 orbifold fibres of $\tau_T$ indeed holds and has been proved by Lotay
in~\cite{lotay2009}.
        \begin{varthrm}[Theorem~A] \mbox{\rm (compare~\cite{lotay2009})}
Let $M$ be a compact 7-manifold with a smooth one-parameter family of
\mbox{$G_2$-structures} given by closed 3-forms $\phi_T\in\Omega^3_+(M)$,
$T>T_0$, defined by the generalized connected sum construction
in~\S\ref{pencils}. Let $\tau_{T}:M\to S^3$, be a coassociative K3
fibration, with respect to $\phi_T$, defined in Theorem~\ref{approx}.
Suppose that $\phi_T+d\eta_T$, is a smooth family of torsion-free
$G_2$-structures on~$M$, such that
$\|\eta_T\|_{L^p_k}<K_{p,k}e^{-\lambda T}$ for each $p>1$, $k\ge 0$.

Then there exists $T_1$ and for any $T>T_1$ and $0\le s\le 1$ a
$C^{1\alpha}$ vector field $v_{T,s}$ on~$M$, smooth away from the singular
fibres of~$\tau_T$ and satisfying
$\|v_{T,s}\|_{L^p_k}<K_{p,k}\,s\,e^{-\lambda T}$
and $\|v_{T,s}\|_{C^1}<K s\,e^{-\lambda T}$,
with support of~$v$ contained in a neighbourhood $U$ of the singular fibres
of~$\tau_T$, and such that $\phi+s\,d\eta_T$ vanishes on every singular
fibre of the perturbed fibration $\tau_T\circ\exp^{-1}_{v_{T,s}}:M\to S^3$.
The neighbourhood $U$ may be chosen not to meet the neck of~$M$, i.e.\
$U\subset (W_1(0)\times S^1\sqcup W_2(0)\times S^1)$.
        \end{varthrm}
We explain in~\S\ref{linear} a framework of appropriate weighted Sobolev
spaces making $(D_1F)_0$ into a Fredholm map
and in~\S\ref{proof} show how the surjectivity of $(D_1F)_0$ then
follows from Lotay's {\em stability} result, in the case of conical
singularities of K3 orbifold fibres. The surjectivity implies that the
deformations of coassociative K3 orbifold fibres of $\tau_T$ are
{\em unobstructed}.

The hypothesis of the next result assumes the assertion of Theorem~A.
We may now deal with an approximating fibration whose singular fibres are
precisely coassociative with respect to a $G_2$-structure $\phi_T$, for
each~$T$.
        \begin{varthrm}[Theorem~B]
Let $M$ be a compact 7-manifold with a smooth one-parameter family of
\mbox{$G_2$-structures} given by closed 3-forms $\phi_T\in\Omega^3_+(M)$,
$T>T_1$, defined by the generalized connected sum construction
in~\S\ref{pencils}. Let $\tau_{T}:M\to S^3$ for $T>T_1$ be a coassociative K3
fibration map, with respect to $\phi_T$, defined in Theorem~\ref{approx}. 

Then there exists $\eps>0$ so that if $T>T_1$ and an
exact form $d\eta\in\Omega^3(M)$ vanishes on every singular fibre
of~$\tau_{T}$ and $\|d\eta\|_{C^0(M)}<\eps$, relative to the metric
$g(\phi_T)$, then there is a unique smooth vector field
$\tilde{v}_{T}(\eta)$ on~$M$ such that:

(i) $\tilde{v}_{T}$ vanishes on the singular fibres of $\tau_{T}$
and is point-wise orthogonal to each smooth fibre $X$ of $\tau_{T}$
and $\tilde{v}_{T}\corner\phi_T|_X$ is $L^2$-orthogonal to the harmonic
self-dual forms on~$X$ relative to the metric $g(\phi_T)|_X$;

(ii) $\tilde{v}_{T}(\eta)$ depends smoothly on $T$ and $d\eta$ and
$\|\tilde{v}_{T}\|_{C^1}=O(\|d\eta\|_{C^1})$;

(iii) $\phi_T+d\eta$ vanishes on the fibres of
$\tau_{T}\circ\exp_{\tilde{v}_T(\eta)}^{-1}$.
        \end{varthrm}
Theorem~B is proved in~\S\ref{thrm.B}.

The estimates of the vector fields $v_{T,s}$ in Theorem~A ensure that
$\exp_{v_{T,s}}$ for any large~$T$ is a well-defined diffeomorphism of~$M$
isotopic to $\id_M$.
Denote $\tilde\phi_T=\exp_{v_{T,s}}^*(\phi_T+d\eta_T)$; then, for $T>T_1$,
the form $\tilde\phi_T$ vanishes on the singular fibres of~$\tau_T$.
The form $\tilde\phi_T$ is in the cohomology class of $\phi_T$ and
$d\tilde\eta_T=\tilde\phi_T-\phi_T$ tends to zero in $C^\infty$ as
$T\to\infty$. For any large~$T$, Theorem~B applies to $d\tilde\eta_T$
and gives a second vector field $\tilde{v}(\tilde{\eta}_T)$ on~$M$,
so that the diffeomorphism $\exp_{\tilde{v}(\tilde{\eta}_T)}$ of~$M$ is
well-defined and fixes the singular fibres of~$\tau_T$.
Then $\phi_T+d\eta_T$ vanishes on the fibres of $\tau_T\circ h_T$, where
$$
h_{T}=\exp_{v_{T,1}}^{-1}\circ\exp_{\tilde{v}_{T,s}}^{-1}.
$$
Thus the fibres of $\tau_T\circ h_T$ are coassociative calibrated by the
4-form $*_{\phi_T+d\eta_T}\phi_T+d\eta_T$ and are minimal submanifolds
of the holonomy-$G_2$ manifold $(M,\phi_T+d\eta_T)$ as required in the Main
Theorem.

\section{Linear analysis on coassociative K3 orbifolds}
\label{linear}

Before going to prove Theorem~A we need to deal with the analytic issues
arising in the deformation theory of the singular, K3 orbifold
coassociative fibres of the map~\eqref{solid.torus}. We begin, in this and
the next section, by showing that the infinitesimal coassociative
deformations of these K3 orbifold fibres are unobstructed.

Our treatment is similar in spirit to one previously used by Joyce in the
series of papers on special Lagrangian submanifolds with conical
singularities, including~\cite{joyce-cone-1,joyce-cone-3}, in that we apply
elliptic theory on non-compact manifolds using weighted Sobolev spaces
and extend these by certain finite-dimensional spaces to eliminate
the obstruction space. Lotay~\cite{lotay} applied the method
of~\cite{joyce-cone-1,joyce-cone-3} and other papers in the same series
to show that deformations of coassociative submanifolds with conical
singularities may in general be obstructed.
Recently Lotay developed a rather general deformation and stability theory
for coassociative submanifolds whose conical singularities arise from
complex cones in~$\CX^3$.

The ordinary double point singularities of complex surfaces in a
threefold~$W$ are an instance of conical singularities. On the other hand,
complex surfaces in $W$ define coassociative submanifolds in $W\times S^1$
with respect to the product $G_2$-structure corresponding to a Calabi--Yau
structure on~$W$. As we explain below, using the stability result
of~\cite{lotay2009}, in the case of ordinary double point singularities
it is possible to set up an unobstructed theory.

Throughout this section, we work on a neighbourhood
$\Mloc=\Wloc\times S^1\subset M$
of a singular fibre, $X_0$ say, of the map~$\tau=\tau_T$ defined
in~\eqref{solid.torus}. Here $\Wloc=(\tau^{(j)})^{-1}(U_0)\subset W_j$
$j=1\text{ or }2$, and $U_0\subset\CP^1$ is an open disc, such that
$X_0\subset\Wloc$ and $\Wloc$ contains no other singular fibres
of~$\tau^{(j)}$. Respectively, $\Mloc=\tau^{-1}(U\times S^1)$;
mark a point $0\in S^1$ and identify $\Wloc$ with
$\Wloc\times\{0\}\subset\Mloc$.
(Note that writing $\tau$ rather $\tau_T$ is justified here as
there are no singular fibres on the neck of~$M$ and the restriction of
the fibration map to a neighbourhood away from the neck does not depend
on~$T$.)
We consider the 7-manifold $\Mloc$ with the product $G_2$-structure
$$
\phi_{CY}=\omega\wedge d\theta + \im\Omega
$$
induced by a Calabi--Yau structure $(\omega,\Omega)$ on~$W$ (as before,
$\theta$ is an `angle coordinate' on~$S^1$), so $\phiCY$ vanishes on the
smooth part of $X_0$. We shall set up a technical framework to deal 
with the local coassociative deformations of $X_0$, extending McLean's
approach.

Recall that $X_0$ is an orbifold degeneration of K3 surface with one
ordinary double point and no other singularities.
Denote by $w_0$ the ordinary double point point of~$X_0$ and by
$X'_0=X_0\minus\{w_0\}$ the complement smooth non-compact complex surface.
We shall always use on a neighbourhood of~$w_0$ in~$W$ the `Morse' local
complex coordinates $(z_1,z_2,z_3)$ discussed in~\S\ref{pencils};
recall that the local expression of $\tau^{(j)}$ in these coordinates is
$z_1^2+z_2^2+z_3^2$ and a neighbourhood of $w_0$ in $X_0$ corresponds
to a neighbourhood of~0 in the cone on~$\RE P^3$
        \begin{equation}\label{odp}
C_0=\{z_1^2+z_2^2+z_3^2=0\}\subset\CX^3.
        \end{equation}
We shall sometimes use on a neighbourhood of $w_0$ in~$X_0$ the
uniformizing coordinates $(u_1,u_2)\in\CX^2$ defined in~\S\ref{pencils},
so that
        \begin{equation}\label{embed}
(z_1,z_2,z_3)=\bigl(i(u_1^2+u_2^2),\; u_1^2-u_2^2,\; 2u_1u_2\bigr).
        \end{equation}

The Ricci-flat K\"ahler metric $\omega$ on~$\Wloc$ defines, by restriction,
an incomplete K\"ahler metric on~$X'_0$ and which can be written near~$w_0$
in the form
        \begin{equation}\label{orbi}
dr^2+r^2g_3+O(r^3),\qquad\text{as }r\to 0.
        \end{equation}
where $r$ is the geodesic polar radius at~$w_0$ in $W$, $g_3$ is some smooth
metric on $\RE P^3$, and $o(r^3)$ is understood in the sense of the uniform
convergence on $\RP^3$ with all derivatives.

\begin{remark}
More explicitly, the metric $g_3$ can be computed by restricting the inner
product on the real tangent space $T_{w_0}W$ induced by $\omega$ to the
intersection of the unit sphere and the tangent cone of $X_0$ at $w_0$.
In the uniformizing coordinates $u_1,u_23$, it is obtained by substituting
the expressions~\eqref{embed} into (the real part of)
$h_{i\bar{j}}dz_id\bar{z}_j$, where $h_{i\bar{j}}$ is the Hermitian inner
product defined by $\omega$ at $w_0$. Note that the K\"ahler form $\omega$
need not be `compatible' with $z_j$'s in any special way.
In particular, $g_3$ is need not in general be induced in the standard way
from the round metric on $S^3$, nor the `obvious' metric on the link of
$C_0$ induced by the Euclidean metric on $\CX^3$.
\end{remark}

Note also that the expression for the metric~\eqref{orbi} in the uniformizing
coordinates degenerates at~$w_0$. In particular, the local 2-form
$du_1\we du_2$ on~$X_0$ is smooth in the orbifold sense but its point-wise
norm relative to the metric $g(\phiCY)$ blows up at~$w_0$,
$|du_1\we du_2|=O(r^{-1})$ as $r\to 0$.

In what follows, we extend the local coordinate $r$
to a positive smooth function, still denoted by~$r$, defined on all of
$X'_0$ and such that $r>1$ away from a coordinate neighbourhood of~$w_0$.

A {\em local deformation of~$X_0$ in~$\Mloc$} is defined as $\exp_v(X_0)$,
where $v$ is a vector field on a neighbourhood $U$ of $X_0$ in~$\Mloc$ with
a small $C^1$-norm, so that $\exp_v:U\to\Mloc$ is a
diffeomorphism of $U$ onto its image.  We regard two local deformations
$\exp_v(X_0)$ and $\exp_{v'}(X_0)$ as equivalent if
$\exp_{v'}(X_0)|_{X_0}=\exp_{v'}(X_0)|_{X_0}\circ \Phi_0$, for some
diffeomorphism $\Phi_0$ of $X_0$ onto itself.  (Note that any diffeomorphism
of the orbifold $X_0$ necessarily fixes~$w_0$.)

The follows is a direct corollary of the tubular neighbourhood
Theorem~\ref{tubular}.
        \begin{prop}\label{tub}
Let $B_x(\rho)$ denote a ball of radius $\rho$ about zero in the
fibre of the normal bundle $N_{X'_0/\Mloc}$ over $x$ with the inner
product induced by the metric on~$\Mloc$.
There exists $\eps>0$ such that the Riemannian exponential map an defines a
diffeomorphism of an open neighbourhood
$U_\eps=\cup_{x\in X'_0}B_x(\eps r)\subset N_{X'_0/\Mloc}$
onto a neighbourhood of $X'_0$ in $\Mloc$, where $r$ is
the `polar radius-function' on $X_0$ defined above.
        \end{prop}

It follows that the bundle isometry $N_{X'_0/\Mloc}\cong\Lambda^+T^*X'_0$
(cf.~\eqref{normal}) bijectively identifies local deformations of $X_0$
fixing~$w_0$ defined by vector fields $v$ point-wise orthogonal to $X_0$ and
the forms $\psi\in\Omega^+(X'_0)$ if the
uniform norms of $r^{-1}v$ and $r^{-1}\psi$ are less than~$\eps$ given
by Proposition~\ref{tub}.
For any such $\psi=v\corner\phiCY|_{X'_0}$, an orbifold~$\exp_v(X_0)$ will be
coassociative if and only if $\omega$ is a zero of `McLean's map'
        \begin{equation}\label{mclean-orbi}
F:\psi=v\lrcorner\phi_{CY}|_{X'_0}\in\Omega^+(X'_0)\to
\exp_v^* \phi_{CY}|_{X'_0}\in\Omega^3(X'_0).
        \end{equation}
Recall from Proposition~\ref{CA} that the linearization of~$F$ at
$\psi=0$ is an overdetermined-elliptic differential operator
$(dF)_0=d:\Omega^+(X_0)\to\Omega^3(X_0)$
and that the image of $F$ consists of exact 3-forms on~$X'_0$.
In order to apply the implicit function theorem to~$F$ we require a
choice of Banach space completions for the space of self-dual forms
$v\lrcorner\phi_{CY}$ on $X'_0$ arising from local deformations of
$X_0$ and also for the space of bounded exact 3-forms on $X'_0$, so that the
exterior derivative extends to a {\em surjective} operator between the
two Banach spaces.

The K\"ahler metric induced on $X'_0$ from $\Wloc$ does {\em not} extend
to a smooth orbifold metric on~$X_0$, so trying to work with orbifold
versions of e.g.\ the usual Sobolev spaces on a compact $X_0$ and is not a
very promising way.
Instead, we use a `conformal blow-up' of~$X'_0$ at the singular point
and apply the elliptic theory for manifolds with asymptotically
cylindrical ends from~\cite{LM,MP,tapsit}.

The non-compact submanifold $X'_0$ is diffeomorphic to a smooth
4-manifold with cylindrical end $\RE_{>0}\times\RP^3$ which corresponds to
a neighbourhood of $w_0$ via $r=e^{-t}$, where $t$ is the coordinate on the
$\RE_{>0}$ factor. Restricting to the end of $X'_0$ we can write
        \begin{equation}\label{conf}
dr^2+r^2g_3+o(r^3)=e^{-2t}\gcyl=
e^{-2t}(dt^2+g_3+o(e^{-t})),\qquad\text{as }t\to\infty,
        \end{equation}
which shows that the metric~\eqref{orbi}
induced on $X'_0$ from~$\Mloc$ is conformally equivalent to an
{\em asymptotically cylindrical} metric $\gcyl$. In particular, the 
self-dual forms defined by the metrics $\gcyl$ and $\omega|_{X'_0}$ are
the same.

We shall need exponentially weighted Sobolev spaces on~$X'_0$.
By definition, $e^{-\delta t}L^p_k(X'_0)_{\text{cyl}}$ is the space
of functions $e^{-\delta t}f$ such that $f\in L^p_k(X'_0)$ and the norm
is defined by $\|e^{-\delta t}f\|_{e^{-\delta t}L^p_k}$ $=\|f\|_{L^p_k}$.
Here we used the subscript `cyl' to indicate that the $L^p_k$ norm in the
previous sentence is calculated using the metric $\gcyl$. This
will be important when we consider the differential forms on~$X'_0$.

There is a preferred choice of weight $\delta=k-4/p$. Denote the
corresponding weighted spaces by~$W^p_k(X_0)$. In terms of the radial
parameter $r$ on~$X'_0$ the $W^p_k$-norm is expressed as
        \begin{equation*}\label{w}
\|f\|_{W^p_{k}(X_0)}=\|\frac{f}{r^{k}}\|_p+
\|\frac{\nabla f}{r^{k-1}}\|_p+\ldots+
\|\nabla^k f\|_p .
        \end{equation*}
The above expression for $W^p_k$ norm extends to the differential forms
of any degree~$m$ on~$X'_0$ but note that the point-wise norms of the
$m$-forms are rescaled by the `conformal weight' factor $e^{-mt}$ when
passing to the  metric~$\gcyl$. In view of this, we define
$$
W^p_{k}\Omega^m(X_0)=
e^{-(k-4/p+m)t}L^p_k\Omega^m(X'_0)_{\text{cyl}}.
$$
Then the exterior derivative extends to a bounded linear map
$W^p_k\Omega^m(X'_0)\to W^p_{k-1}\Omega^{m+1}(X'_0)$.

There is a simple relation between the $W^p_k$ spaces and the
usual, unweighted Sobolev spaces $L^p_k(X_0)$. The following result is
proved in~\cite{Biquard} in the case of a flat Euclidean ball (with $r$ the
Euclidean distance to the origin). However, the argument of the proof works,
with only a change of notation, for a punctured ball endowed with a metric
having the `conical' form~\eqref{orbi}. Considering a punctured
neighbourhood of $w_0$ in $X_0$ as the quotient of a 4-dimensional
punctured ball with respective $\pm 1$-invariant metric and restricting
attention to $\pm 1$-invariant functions we obtain.
        \begin{prop}[cf.~\mbox{\cite[Theorem 1.3]{Biquard}}]\label{biq}
Suppose that $\ell$ is a non-negative integer such that
$\ell-1<k-4/p<\ell$ and let $\delta=k-4/p$. Then one has
        \begin{equation}\label{unweight}
W^p_{k}(X_0)=\{f\in L^p_k(X_0)\; |\;
\lim_{w\to w_0}\nabla^m f(w)=0\text{, for all }0\le m\le\ell-1\}
        \end{equation}
and the $W^p_k$-norm on the left-hand side is equivalent to the
$L^p_k$-norm on the right-hand side.
        \end{prop}
The vanishing condition in the right-hand side of~\eqref{unweight} makes
sense as $L^p_k$ embeds in~$C^{\ell-1}$. Proposition~\ref{biq} extends in
the usual way to sections of vector bundles over~$X_0$ by considering a
$\pm1$-equivariant local trivialization near~$w_0$.

In view of the local regularity results for coassociative calibrated
manifolds~\cite[\S~IV.2.7]{calibrated} we require a Banach space
consisting of the $C^1$ self-dual forms on~$X'_0$.
We shall use the completion of the space of self-dual
forms in a $W^p_{k}$-norm fixing a choice of 
$p>1,\; k\in\ZE$ such that
        \begin{equation}\label{sobolev}
1<k-4/p<2 
        \end{equation}
With this choice, every self-dual form $\psi\in W^p_k\Omega^+(X'_0)$
vanishes to order two at $w_0$, in particular $\psi$ is Lipschitz
continuous at $w_0$ with any Lipschitz constant $\eps>0$.
Whenever the $W^p_k$ norm of $\psi$ is small, the corresponding section $v$
of $N_{X'_0/\Mloc}$ defines a local deformation $\exp_v$ of~$X_0$ which
fixes the singular point~$w_0$ and the tangent cone at~$w_0$.

In order to include local deformations of~$X_0$ which move the tangent
cone and the singular point, we extend the weighted Sobolev space of
self-dual forms by adding a finite-dimensional space
$$
E_0=\{\vea\lrcorner\phi_{CY}|_{X'_0} :
\mathbf{e}\in T_{w_0}\Mloc,\; L \in\End(T_{w_0}\Mloc) \}. 
$$
Here $\vea$ denotes a choice of a smooth vector field on $\Mloc$
smoothly depending on the parameters~$\mathbf{e},L$,
such that
$\vea(w_0)=\mathbf{e}$ and $d\vea (w_0)= L$.
In the last condition we define $d\vea$ by using the coordinates
$z_i,\theta$ near $w_0\in\Mloc$ to express $\vea$ locally as a smooth map
from a neighbourhood of zero in $T_{w_0}$ to~$T_{w_0}$.
The Banach space $E_0+W^p_{k}\Omega^+(X'_0)$, for $1<k-4/p<2$,
does not depend on the choice of vector fields~$\vea$ given above
as the ambiguity is $O(r^2)$, $r\to 0$, which is contained in~$W^p_k$ by
Proposition~\ref{biq}.
It is not difficult to check, using Proposition~\ref{biq}, that
for every smooth local deformation $\exp_v$ of~$X_0$ the self-dual form
$v\lrcorner\phiCY|_{X'_0}$ is in
$E_0+e^{-\delta t}L^p_{k}\Omega^+(X'_0)_{\text{cyl}}$.
On the other hand, straightforward calculation in local coordinates on
$\Mloc$ near $w_0$ shows that any smooth exact 3-form on $\Mloc$ restricts
to a form in $d(E_0)+W^p_{k-1}(\Omega^3(X'_0)$.
(As $k-1-4/p>0$, the space $W^p_{k-1}(\Omega^3(X'_0)$ alone only contains
forms vanishing at~$w_0$.)

We are now ready to state the main result of this section.
        \begin{thrm}\label{noobstr}
If the Sobolev space parameters satisfy~\eqref{sobolev} then the exterior
derivative defines a bounded linear map between Banach spaces
        \begin{equation}\label{lin.orbi}
d: E_0+W^p_{k}\Omega^+(X'_0)\to
d(E_0)+\{\eta\in W^p_{k-1}\Omega^3(X'_0):d\eta=0\} 
        \end{equation}
which is surjective and has a one-dimensional kernel spanned
by the K\"ahler form $\omega|_{X'_0}$.
        \end{thrm}
        \begin{remark}
The 1-dimensional kernel of~\eqref{lin.orbi} arises from the
$S^1$-symmetry of the torsion free $G_2$-structure $\phi_{CY}$
on $W\times S^1$. 
Thus Theorem~\ref{noobstr} shows that the $S^1$-families of coassociative
K3 orbifolds arising in the approximating fibration~\eqref{solid.torus} are
{\em maximal} deformation families (the smooth fibres, of course, have the
same property by Theorem~\ref{CA} and~\cite{mclean} as $b^+=3$ for a K3
surface).
        \end{remark}

\section{Proof of Theorem~\ref{noobstr}}
\label{proof}

\begin{remark}
Putting $x=e^{-t}$ we can think of the manifold $X'_0$ as the
interior of a compact manifold with boundary
$\overline{X'_0}=X'_0\cup(\{x=0\}\times\RP^3)$
obtained by adding a copy of~$\RP^3=\{x=0\}$ `at infinity'.
(The added $\RP^3$ may also be canonically identified with the unit
spherical space form $S^3/\pm1$ in the tangent cone of $X_0\subset\Wloc$
at~$w_0$.) The asymptotically cylindrical metric $\gcyl$ on~$\overline{X'_0}$
can be written in the form 
        \begin{equation}\label{bmetr}
\gcyl=\frac{dx^2}{x^2}+\tilde{g},
        \end{equation}
where $\tilde{g}$ is a symmetric semi-positive definite form smooth up
to the boundary:  at any point in~$\{x=0\}$, $\tilde{g}$~is smooth in the
$\RP^3$ directions and has one-sided derivatives in~$x$ of any order.
Then~\eqref{bmetr} gives an instance of a `smooth exact $b$-metric', as
defined by Melrose \cite[Ch.~2]{tapsit}. The results proved
in~\cite{tapsit} for manifolds with smooth exact $b$-metrics can therefore
be applied to $(X'_0,\gcyl)$.
\end{remark}

In will be convenient to reduce Theorem~\ref{noobstr} to a
result concerning a linear operator between the $W^p_k$-spaces alone.
We look at the kernel first.
(The subscript `cyl' at the exponentially weighted spaces on~$X'_0$
will now be dropped from the notation.)
        \begin{prop}\label{1-1}
Suppose that $k-4/p>1$. The map~$d:e^{-\delta t}L^p_k\Omega^+(X'_0)\to 
e^{-\delta t}L^p_{k-1}\Omega^3(X'_0)$ has a 3-dimensional kernel for
any $0\le\delta<1$, a 1-dimensional kernel for any $1\le\delta<2$, and is
injective for any $\delta\ge 2$. 
        \end{prop}
        \begin{pf}
This is an application of the Hodge theory on asymptotically cylindrical
manifolds. The kernel of~$d$ is contained in the kernel of the
Laplacian acting on $e^{-\delta t}L^p_k\Omega^+(X'_0)$. Furthermore, the
standard integration by parts argument is valid when $\delta>0$ and shows
that the two kernels coincide. By~\cite[Propn.~4.9]{aps} or
\cite[Propn.~6.14]{tapsit}, the $L^2$-kernel of the Laplacian on the $m$-forms
on an asymptotically cylindrical manifold~$X'_0$ is isomorphic to the
image of the natural inclusion homomorphism $H_c^m(X'_0)\to H^m(X'_0)$
of the de Rham cohomology groups, where the subscript `c' indicates the
cohomology with compact support.

Considering the exact sequence of the de Rham cohomology groups
$$
\ldots\to H^{m-1}(\RP^3)\to H_c^m(X'_0)\to H^m(X'_0)\to H^m(\RP^3)\to\ldots
$$ 
we find that $H^0_c(X'_0)=0$, $H^4_c(X'_0)\cong\RE$ and the inclusion
homomorphism $H^m_c(X'_0)\to H^m(X'_0)$ is an isomorphism for $1\le m\le3$.
Recall from~\S\ref{pencils} that $X'_0$ is isomorphic to the complement of
a $(-2)$-curve in a K3 surface, $X$ say; this $(-2)$-curve is topologically a
sphere. Considering the Maier--Vietoris exact sequence for the union of
$X'_0$ and a tubular neighbourhood of the `missing' $(-2)$-curve we find
that  $H^2(X'_0)$ is a complement in $H^2(X)$ of the one-dimensional
subspace generated by the Poincar\'e dual of the $(-2)$-curve. The
cup-product on $H^2_c(X'_0)$ has maximal positive subspace $H^+(X'_0)$ of
dimension~3 and so the $L^2$-kernel of the Laplacian on $\Omega^+(X'_0)$
is 3-dimensional.

It is easy to identify, for $0<\delta<1$, the 3-dimensional space of
$O(e^{-\delta t})$ exponentially decaying closed self-dual forms
on~$X'_0$. This space is spanned by the restriction $\omega|_{X_0}$ of the
K\"ahler form on~$\Wloc$ and the restrictions of the real and imaginary parts
of the holomorphic $(2,0)$-form
$\bigl((\p\tau)^\sharp \lrcorner\Omega\bigr)|_{X^0}$,
defined using is the K3 fibration map $\tau$ on~$\Wloc$.
The K\"ahler form $\omega|_{X_0}$, measured with the metric $\gcyl$,
is $O(e^{-2t})$, as $t\to\infty$, but not $O(e^{-(2+\eps)t})$ for any
$\eps>0$. As the $(1,0)$-form $\p\tau$ has a zero of order~1 at $w_0$,
the vector field $(\p\tau)^\sharp$ is $O(1/r)$, as $r\to 0$, when
measured in to the K\"ahler metric $\omega|_{X_0}$. Hence the real and
imaginary parts of $\bigl((\p\tau)^\sharp \lrcorner\Omega\bigr)|_{X^0}$
are $O(e^{-t})$ but not $O(e^{-(1+\eps)t})$in the metric~$\gcyl$.

Neither of the three self-dual forms spanning the $L^2$-kernel of the
Laplacian (and the $L^2$-kernel of~$A$) is in $e^{-\delta t}L^p_k$ for
$\delta>2$.
        \end{pf}
        \begin{cor}\label{1-1c}
If $1<k-4/p<2$ then the map~\ref{lin.orbi} has a one-dimensional
kernel spanned by the K\"ahler form~$\omega|_{X'_0}$.
        \end{cor}
        \begin{pf}
It is not difficult to check that
$\omega|_{X'_0}\in E_0+W^p_{k}\Omega^+(X'_0)\subset
e^{-\delta t}L^p_k\Omega^+(X'_0)$,
for any $1<\delta<2$ but $E_0+W^p_{k}\Omega^+(X'_0)$.
        \end{pf}
The next two lemmas determine the codimensions of the relevant
$W^p_k$ spaces.
        \begin{lemma}\label{codim}
If $1<k-4/p<2$ then the codimension of $W^p_k\Omega^+(X'_0)$ in
$E_0+W^p_k\Omega^+(X'_0)$ is~$41$.
        \end{lemma}
        \begin{pf}
If $\mathbf{e}\neq 0$ then $|e^{2t}\vea\corner\phiCY|_{\gcyl}$ has a
non-zero lower bound on an open subset of the cylindrical end,
of the form $\RE_+\times U'$ where $U'$ is open in $\RE^3$. So
$\vea\corner\phiCY$ with $\mathbf{e}\neq 0$ is never in
$W^p_k\Omega^+(X'_0)$ if $k-4/p>0$.

The form 
$v_{0,L}\corner\phiCY$ will be in $W^p_k\Omega^+(X'_0)$ with $1<k-4/p<2$
precisely if $L$ leaves invariant the tangent cone of $X_0$ at~$w_0$. This,
in turn, will be the case if and only if the Zariski tangent space
$T_{w_0}\Wloc$ of $X_0$ is an invariant subspace of~$L$ and the restriction
of~$L$ to $T_{w_0}\Wloc$ is up to a complex factor an element of $SO(3,\CX)$
(so that $L$ preserves the tangent cone~\eqref{odp}).
We find that the subspace of endomorphisms in $\End T_{w_0}\Mloc$
preserving the tangent cone has dimension $15$. As the dimension of~$E_0$
is $7+49$ the result follows.
        \end{pf}
        \begin{lemma}\label{target}
If $1<k-4/p<2$ then the codimension of $W^p_{k-1}\Omega^3(X'_0)\cap\Ker d$
in $d(E_0)+(W^p_{k-1}\Omega^3(X'_0)\cap\Ker d)$ is~$18$.
        \end{lemma}
        \begin{pf}
Recall from \S\ref{linear} that the space
$d(E_0)+(W^p_{k-1}\Omega^3(X'_0)\cap\Ker d)$ contains the restrictions to
$X'_0$ of all the smooth exact 3-forms on~$\Mloc$. On the other hand,
$W^p_{k-1}\Omega^3(X'_0)$ is precisely the space of $L^p_{k-1}$ 3-forms on
$X'_0$ with zero limit at~$w_0$. Since $X_0\subset\Wloc$ and the Zariski
tangent space of~$X_0$ at $w_0$ is $T_{w_0}\Wloc$ the codimension of
interest may be computed as $\dim\Lambda^3T_{w_0}^*\Wloc - 
\dim(d(E_0)\cap W^p_{k-1}\Omega^3(X'_0))$.

A calculation in the local complex coordinates on $\Wloc$ near $w_0$ shows
that\linebreak
$d(\vea\corner\phiCY)|_{X'_0}\in W^p_{k-1}\Omega^3(X'_0)$ holds
precisely if $d(\vea\corner\phiCY)|_{\Wloc}(w_0)$ is 
the real or imaginary part of a $(3,0)$-form on $T_{w_0}\Wloc$.
        \end{pf}
From Corollary~\ref{1-1c} and Lemmas~\ref{codim},~\ref{target} and some
straightforward linear algebra we find that Theorem~\ref{noobstr} is
equivalent to the following technical result on an asymptotically
cylindrical 4-manifold $(X'_0,\gcyl)$.
        \begin{varthrm}[Theorem~\ref{noobstr}$\, '$]
If $1<k-4/p<2$ then the injective linear map
        \begin{equation}\label{A}
A:\psi\in W^p_k\Omega^+(X'_0)\to
d\psi\in W^p_{k-1}\Omega^3(X'_0)\cap\Ker d,
        \end{equation}
has a 22-dimensional cokernel.
        \end{varthrm}
In the remainder of this section we prove Theorem~\ref{noobstr}$\, '$.
Note that the injectivity of $A$ follows from Proposition~\ref{1-1}.
Therefore, the dimension of $\Coker A$ is minus the Fredholm index of~$A$.
For the index computation, it is convenient to observe that the map $A$ is
equivalent to a component of an elliptic operator
        \begin{equation}\label{YM}
D: (f,\psi)\in e^{-\delta t}L^p_k(\Omega^0\oplus\Omega^+)(X'_0)\to
df-*d\psi\in e^{-\delta t}L^p_{k-1}\Omega^1(X'_0).
        \end{equation}
For any real $\delta$, we write $\Ind_{-\delta} D$ to indicate that
$D$ is considered on the $e^{-\delta t}$-weighted Sobolev spaces. A
similar notation will be used for $A$ and for the kernels and cokernels.
        \begin{prop}\label{index}
$\Ind_{-\delta} A = \Ind_{-\delta} D$ for any $\delta>0$.
        \end{prop}
Proposition~\ref{index} will be deduced from the following.
        \begin{lemma}\label{onto}
If $\eps>0$ is sufficiently small then $\Coker_{-\eps}D=\{0\}$.
        \end{lemma}
        \begin{pf}[Proof of Proposition~\ref{index} assuming Lemma~\ref{onto}]
For any $\delta>0$, the cohomology of the weighted de Rham complex on $X'_0$
        \begin{equation}\label{dR}
e^{-\delta t}L^p_k\Omega^0(X'_0)\overset{d}\to 
e^{-\delta t}L^p_{k-1}\Omega^1(X'_0)\overset{d}\to
e^{-\delta t}L^p_{k-2}\Omega^2(X'_0)\overset{d}\to\ldots
        \end{equation}
is isomorphic to the de Rham cohomology with compact
support $H^*_c(X'_0)$  \cite[Propn.~6.13]{tapsit}.
The formal $L^2$-adjoint of~\eqref{dR} is a complex of the
$e^{\delta t}$-weighted spaces. The cohomology of the latter complex at the
$\Omega^m$-term is isomorphic to the de Rham cohomology~$H^{4-m}(X'_0)$.

Recall from the proof of Proposition~\ref{1-1} that
$H^0_c(X'_0)=H^1_c(X'_0)=0$. So the exterior derivative maps
$e^{-\delta t}L^p_k\Omega^0(X'_0)$, for each $\delta>0$, isomorphically
onto the subspace of closed forms in
$e^{-\delta t}L^p_{k-1}\Omega^1(X'_0)$.
Standard integration by parts argument of the Hodge
theory for compact manifolds is valid for the exponentially decaying forms
on $X'_0$ and shows that the spaces of 1-forms
$d(e^{-\delta t}L^p_{k}\Omega^0(X'_0))$ and
$d^*(e^{-\delta t}L^p_{k}\Omega^+(X'_0))$ are $L^2$-orthogonal if $\delta>0$.
By Lemma~\ref{onto}, $D$ is a surjective operator between
$e^{-\eps t}$-weighted spaces for any small $\eps>0$, so we obtain a
decomposition
$$
e^{-\eps t}L^p_{k-1}\Omega^1(X'_0)
= d(e^{-\eps t}L^p_k\Omega^0(X'_0))
\oplus d*(e^{-\eps t}L^p_k\Omega^+(X'_0)).
$$
For an arbitrary $\delta>0$ we can write any
$\xi\in e^{-\delta t}L^p_{k-1}\Omega^1(X'_0)$ as $\xi=d\xi_0+*d\xi_+$
where $\xi_+\in e^{-\eps t}L^p_k\Omega^+(X'_0)$ for some {\em small}
$\eps>0$ but $\xi_0\in e^{-\delta t}L^p_k\Omega^0(X'_0)$
(as explained above).
Therefore, the image of
$D$ on the $e^{-\delta t}$-weighted spaces is complemented by
$d^*$-closed forms and the proposition follows.
        \end{pf}
For the proof of Lemma~\ref{onto}, we need to recall from \cite{LM}
or~\cite{tapsit} some Fredholm theory for elliptic operators on
asymptotically cylindrical Riemannian manifolds.
Let $\Lambda_\infty^+ T^*(\RE_{>0}\times\RP^3)$ denote the bundle of
self-dual forms with respect to the metric $dt^3+g_3$, the asymptotic model
of~$\gcyl$. A point-wise orthogonal projection, relative to $dt^3+g_3$,
$\sigma:\Lambda^+ T^*(\RE_{>0}\times\RP^3)\to
\Lambda_\infty^+ T^*(\RE_{>0}\times\RP^3)$ defines a bundle isomorphism
asymptotic to the identity as $t\to\infty$.
The coefficients of $D$ are determined by the metric $\gcyl$ and $D$ is
asymptotic, on the end of~$X'_0$, to an operator
$D_\infty\circ(1\oplus\sigma)$,
where $D_\infty$ is given by the same formula as in~\eqref{YM} but using
the product cylindrical metric $dt^3+g_3$ rather than~$\gcyl$. 
The coefficients of the operator $D_\infty$ on $\RE\times\RP^3$ are
independent of~$t$.
        \begin{prop}\label{b-diff}
(i) The elliptic operator $D$ between $e^{-\lambda t}$-weighted Sobolev
spaces, is Fredholm if and only if $d(\lambda)=0$, where
        \begin{equation}\label{spb}
d(\lambda)=\dim\{e^{-\lambda t}p(t,y)\; |\; p(t,y)\text{ is polynomial in }t
\text{ and }D_\infty(e^{-\lambda t}p(y,t))=0\}.
        \end{equation}
The set $\{\delta\in\RE:d(\delta)\neq 0\}$ is discrete in~$\RE$.

(ii)
The index of $D$ is independent of~$p,k$ but depends on the
weight parameter $\lambda$ according to the formula
        \begin{equation}\label{jump}
\Ind_{\delta''} D-\Ind_{\delta'} D =
\sum_{\delta'<\lambda<\delta''}d(\lambda)
        \end{equation}
for any $\delta'<\delta''$ such that $d(\delta')\neq 0$, $d(\delta'')\neq 0$.

(iii) The kernel of~$D$ consists of smooth forms and is independent of~$p,k$.
The cokernel of~$D$ can be identified with the kernel of the formal
$L^2$-adjoint of~$D$ with respect to~$\gcyl$,
$$
D^*=d^*\oplus d^+: e^{\delta t}L^p_{k-1}\Omega^1(X'_0)\to
e^{\delta t}L^p_k(\Omega^0\oplus\Omega^+)(X'_0).
$$
In particular,
$\dim\Coker_{-\lambda}D=\dim\Ker_{\lambda}D^*$ if $d(\lambda)=0$.
        \end{prop}
Proposition~\ref{b-diff} is a direct application of \cite{LM}
or~\cite{tapsit}.  In the case of~$D$ we don't need to worry 
about the possibility of $\lambda\in\CX\minus\RE$ with
$d(\lambda)\neq 0$, as we shall see in a moment.
        \begin{pf}[Proof of Lemma~\ref{onto}]
The operator $D_\infty$ can be expressed, using some vector bundle
isomorphisms, as $\p_t+D^{(3)}$ where $D^{(3)}$ is a formally
self-adjoint operator on $\RP^3$ whose square is the Laplacian on
$(\Omega^0\oplus\Omega^1)(\RP^3)$ (cf.~e.g.~\cite[pp.~132--134]{mmr}).
This yields $d(0)=1$.
Therefore, $D$ is a Fredholm operator on $e^{-\eps t}$- and on
$e^{\eps t}$-weighted spaces whenever $\eps>0$ is sufficiently small.
Furthermore, for a small $\eps>0$, we obtain
$$
\dim\Ker_{-\eps} D = 3,\quad \dim\Ker_{-\eps} D^*=0
$$
by the Hodge theory arguments as in the proof of Proposition~\ref{1-1},
using also the relation $DD^*=dd^*+\pol d^*d$. On the other hand,
$$
\Ind_{\eps} D - \Ind_{-\eps} D = 1,
$$
by~\eqref{jump}, whence
$$
\dim_{-\eps}\Coker D + \dim_{\eps}\Ker D = 4.
$$
by Proposition~\ref{b-diff}(iii).
But $\dim_{\eps}\Ker D\ge 1+\dim_{-\eps}\Ker D=4$ as the
$e^{\eps t}L^p_k$-kernel contains the constants,
so $\Coker_{-\eps} D=\{0\}$.
        \end{pf}

The next result is proved by Lotay~\cite{lotay2009} and is essentially the
stability property for the coassociative cones defined by complex cones
biholomorphic to~$C_0$.
Clause (i) will be needed for the $C^{1,\alpha}$-regularity claim in
the next section. 

        \begin{prop}[J. D. Lotay]\label{dimensions}
Let $D$ be the elliptic operator defined in~\eqref{YM} over a K3
orbifold $X_0$ and let $d(\lambda)$ be as defined in~\eqref{spb}. Then

(i) any $0<\lambda\le 0$ such that $d(\lambda)\neq 0$ is an integer;

(ii) $\sum_{0<\lambda<\delta} d(\lambda) = 25$, for some $3<\delta<4$.
        \end{prop}

Theorem~\ref{noobstr}$\, '$  follows from Proposition~\ref{dimensions} by
dimension counting. For if $3<\delta<4$ then using ~\eqref{jump},
Propositions~\ref{1-1} and~\ref{index}, and
Proposition~\ref{dimensions}(ii) we obtain $d(\delta)\neq 0$
and
$$
\dim\Coker_{-\delta} A = -\Ind_{-\delta} A = - \Ind_{-\delta} D
= \sum_{0<\lambda<\delta} d(\lambda) - 3 = 22
$$
as required.

{\it Remarks on Proposition~\ref{dimensions}.}
The conformal factor at the Hodge star in~\eqref{YM} is precisely the
inverse of the conformal weight of 1-forms. Therefore, the
operator~$D$ and its asymptotic model $D_\infty$ are
{\em conformally invariant} and can be interchangeably considered on $X'_0$
with the K\"ahler metric $\omega|_{X'_0}$.

The asymptotic model for the cylindrical end of $X'_0$ corresponds
to the cone $C_0$ in
{$T_{w_0}\Wloc\cong\CX^3$} defined by~\eqref{odp} and the asymptotic
model for $\gcyl$ is conformally equivalent to the
K\"ahler metric $dr^2+r^2g_3$ on $C_0$ induced by restricting
from~$T_{w_0}\Wloc$ a Hermitian inner product defined by $\omega(w_0)$.
The kernel elements of $D_\infty$ contributing to the
dimensions $d(\lambda)$ are expressed on $C_0$ as pairs
        \begin{equation}\label{degr}
(r^{\lambda}p_0(\log r, y),r^{\lambda-2}p_+(\log r, y)
\in(\Omega^0\oplus\Omega^+)(C_0\minus\{0\}),
        \end{equation}
where a 0-form $p_0$ and a self-dual form $p_+$ are
polynomial in $\log r$ and smooth in $y\in\RP^3$.
As noted earlier, the conformal rescaling produces an extra factor
$r^{-2}$ for the self-dual forms in~\eqref{degr}.

For a K\"ahler surface $C_0\minus\{0\}$, there is a commutative diagram,
cf.~\cite[\S 4.3.3]{FM},
        \begin{equation}\label{surf}
        \begin{CD}
\Omega^0\oplus\Omega^+      @>{D_\infty}>>          \Omega^1\\
@V{\iota,j+\pi^{0,2}}V{\cong}V                      @V{\cong}V{\pi^{0,1}}V\\
\Omega^0\oplus\Omega^{0,2}  @>{2\bar\p+\bar\p^*}>>  \Omega^{0,1},
        \end{CD}
        \end{equation}
where $\Omega^0$ in the bottom row denotes the 
complex-valued 0-forms, $\iota$ is the inclusion map,
$j(\alpha)=-i\alpha.\omega_{C_0}$ is the contraction with the K\"ahler form
and the inverses of $\pi^{p,q}$ map the respective complex vector space to its
underlying real vector space.

Suppose that
\begin{equation}\label{f.psi}
2\bar\p f+\bar\p^*\psi=0
\end{equation}
and that
$(f,\psi)\in\Omega^0\oplus\Omega^{0,2}$ is obtained from~\eqref{degr}
for some $\lambda>0$ via the left column of~\eqref{surf}. Then since
$\re f$ has degree $r^\lambda$, whereas $\im f$ and
$\psi$ both have degree $r^{\lambda-2}$ we deduce that $\re f=0$. 

If $(f,\psi)$ is a solution of~\eqref{f.psi} then, writing
$\psi=a du_1\we du_2$, applying, respectively, $\p$ and $\p^*$ to the
equation we deduce that $\bar\p^*\bar\p f=0$ and $\p^*\p a=0$,
that is, $f$ and $a$ are harmonic functions on $C_0$.
Therefore, $f$ and $a$ each factorize as $r^\lambda G(y)$ (we omit the
details but cf.~\cite[Propn.~2.4]{joyce-cone-1}).
A function $r^{\lambda-2}G(y)$ is harmonic if and only if $G$ if an
eigenfunction of the Laplacian on the link $\RP^3=\{r=1\}$ of the
singularity of~$C_0$. The eigenvalue of $G$ is then $\lambda(\lambda-2)$
and must be non-negative, therefore there are no solutions $f,\psi$
to~\eqref{f.psi} of the form \eqref{degr} with $0<\lambda<2$.

The harmonic $(0,2)$-form~$\psi$
may be written in the uniformizing coordinates as
$\psi=a\,d\bar{u_1}\we d\bar{u_2}$, for some harmonic complex function~$a$.
Note that $d\bar{u}_1\we d\bar{u}_2$ is homogeneous of degree $-1$ in $r$
and so $a$ factorizes as $r^{\lambda-1} G(y)$.

If we assume that both $f$ and $a$ are {\em homogeneous polynomials}
of even degree in $u_1,u_2,\bar{u}_1$,$\bar{u}_2$ then solving \eqref{f.psi}
with $2\le\lambda<4$ in this case becomes an elementary calculation.
We obtain that $\im f$ must be a constant if $\lambda=2$ and
a harmonic homogeneous quadratic polynomial if $\lambda=3$. This
contributes one real dimension to $d(2)$ and 6 real dimensions to $d(3)$ as
a polynomial $h$ can be found to satisfy \eqref{f.psi}. But $h$ then
is determined up to an anti-holomorphic homogeneous polynomial of
$\bar{u}_1$,$\bar{u}_2$ of degree $2\lambda-2$. 
This contributes further $4\lambda-2$ real dimensions to $d(\lambda)$.
We thus find that $d(1)\ge 2$, $d(2)\ge 7$, and $d(3)\ge 16$. 
Proposition~\ref{dimensions}, in this context, asserts that these are
{\em all} the solutions of~\eqref{f.psi} for $0<\lambda\le 3$ and the
above estimates for $d(\lambda)$'s are in fact equalities.

\section{Perturbing the singular fibres}
\label{thrm.A}

We can now deduce Theorem~A.
        \begin{varthrm}[Theorem~A]
Let $M$ be a compact 7-manifold with a smooth one-parameter family of
\mbox{$G_2$-structures} given by closed 3-forms $\phi_T\in\Omega^3_+(M)$,
$T>T_0$, defined by the generalized connected sum construction
in~\S\ref{pencils}. Let $\tau_{T}:M\to S^3$, be a coassociative K3
fibration, with respect to $\phi_T$, defined in Theorem~\ref{approx}.
Suppose that $\phi_T+d\eta_T$, is a smooth family of torsion-free
$G_2$-structures on~$M$, such that
$\|\eta_T\|_{L^p_k}<K_{p,k}e^{-\lambda T}$ for each $p>1$, $k\ge 0$.

Then there exists $T_1$ and for any $T>T_1$ and $0\le s\le 1$ a
$C^{1\alpha}$ vector field $v_{T,s}$ on~$M$, smooth away from the singular
fibres of~$\tau_T$ and satisfying
$\|v_{T,s}\|_{L^p_k}<K_{p,k}\,s\,e^{-\lambda T}$
and $\|v_{T,s}\|_{C^1}<K s\,e^{-\lambda T}$,
with support of~$v$ contained in a neighbourhood $U$ of the singular fibres
of~$\tau_T$, and such that $\phi+s\,d\eta_T$ vanishes on every singular
fibre of the perturbed fibration $\tau_T\circ\exp^{-1}_{v_{T,s}}:M\to S^3$.
The neighbourhood $U$ may be chosen not to meet the neck of~$M$, i.e.\
$U\subset (W_1(0)\times S^1\sqcup W_2(0)\times S^1)$.
        \end{varthrm}
Theorem~A will be deduced by building up on the
analysis of the previous two sections. The orbifold fibres of $\tau_T$ are
parameterized by finitely many disjoint circles. It will suffice to
restrict attention to one such $S^1$-family. We assume the notation
$\Mloc=\Wloc\times S^1$, $\phiCY\in\Omega^3_+(\Mloc)$, $\tau=\tau_T$,
$X_0\subset\Wloc$, $X'_0=X_0\minus\{w_0\}$, and the radial coordinate
$r$ on $X'_0$, as in~\S\ref{linear}.
For each $\theta\in S^1$, define a finite-dimensional space of smooth
vector fields
$$
\tilde{E}_\theta=
\{\vea|_{X'_0\times\{\theta\}}:
\mathbf{e}\in T_{(w_0,\theta)}(\Wloc\times\{\theta\})\},
$$
with $\vea$ as defined in Theorem~\ref{noobstr}.
        \begin{thrm}\label{extn}
Suppose that $1<k-4/p<2$.
Then there exists $\eps_0>0$ such that for each
$\eta\in\Omega^2(\Mloc)$ with $\|d\eta\|_{C^1(\Mloc)}<\eps_0$ and
$\theta\in S^1$ there is a unique vector field along $X'_0\times\{\theta\}$,
$v(\eta,\theta)\in\tilde{E}_\theta+W^p_{k}N_{(X'_0\times\{\theta\})/\Mloc}$,
with the following properties:

(i) $v(0,\theta)=0$ and $v$ depends smoothly on $\eta$ and $\theta$,

(ii) $\phi_{CY}+d\eta$ vanishes on 
$\exp_{v(\eta,\theta)}(X'_0\times\{\theta\})$.

(iii) the vector field $v(\eta,\cdot)$ on $X'_0\times S^1$ is the
restriction of some $C^{1,\alpha}$ vector field $v(\eta)$ on~$\Mloc$ smooth
away from $X_0\times S^1$,
such that $\|v(\eta)\|_{L^p_k}<C_{p,k}\|d\eta\|_{L^p_k}$
        \end{thrm}
Theorem A is a consequence of Theorem~\ref{extn}.
Choose $T_1$ so that the $C^1$ norm of $d\eta_T$, for $T>T_1$, is smaller
than $\eps_0$ given by Theorem~\ref{extn} for each $S^1$-family of the
singular fibres of~$\tau_T$.  Further increasing $T_1$ if necessary, we can
ensure that the vector field $v(\eta)$ can be chosen with compact
support contained in~$\Mloc$ and with a small $C^1$-norm so that
$\exp_{v(\eta)}$ is a diffeomorphism. Let $U$ be the union of the
neighbourhoods $\Mloc$ for all $S^1$-families of the singular fibres of
$\tau$. The vector field $v_{T,s}$ required in Theorem~A is then obtained
by putting $\eta=s\,d\eta_T$ and extending the $v(\eta)$
by zero on $M\minus U$.

In the rest of this section we prove Theorem~\ref{extn}.

The proof is based on an application of the Implicit Function
Theorem in Banach spaces 
to a `parametric version' of McLean's map
for the family $X_0\times S^1$ of coassociative K3 orbifolds
        \begin{multline}\label{F}
F:\; (v,\eta,\theta)\in
(\tilde{E}_0+W^p_{k}N_{X'_0/\Mloc})
\times\Omega^2(\Mloc)\times S^1
\to\\
\exp^*_{v}(\phi_{CY}+R_\theta^* d\eta)|_{X'_0}
\in d(\tilde{E}_0\corner\phiCY)+(W^p_{k-1}\Omega^3(X'_0)\cap\Ker d),
        \end{multline}
where $R_\theta:\Wloc\times S^1\to \Wloc\times S^1$ denotes an isometry
given by the translations by $\theta$ on the $S^1$ factor (here we
exploited the $S^1$-symmetry to make the Banach space of self-dual forms
independent of~$\theta$). A choice of Banach space for $\eta$ is not
particularly important here as long as it controls e.g.\ the $C^{3}$ norm.
We need the following.
        \begin{prop}\label{smooth}
The map $F$ defined in~\ref{F} is a smooth map of Banach spaces
at any $(v,\eta,\theta)$ with sufficiently small $\|v\|$ in
$\tilde{E}_0+W^p_{k}N_{X'_0/\Mloc}$.
        \end{prop}
        \begin{pf}
The map $F$ is linear in $\eta$ and has an equivariant property
        \begin{equation}\label{equivar}
F(v,\eta,\theta)=F(v,R_\theta^*\eta,0).
        \end{equation}
Therefore, $F$ is
smooth in $\eta,\theta$ and it remains to show that $F$ is smooth in~$v$. 

A related argument was carried out by Baier
in~\cite[Theorem~2.2.15]{baier}, see also~\cite[Theorem~2.5]{joyce-salur}.
We may disregard a finite-dimensional component $\tilde{E}_0$ and pretend
that $v\in W^p_{k}N_{X'_0/\Mloc}\subset C^1 N_{X'_0/\Mloc}$, so $v$ has a
point-wise estimate $|v|<\eps r$ near $w_0$ for any $\eps>0$ (in the metric
$g(\phiCY)$). By Proposition~\ref{tub}, the Riemannian exponential
map is a diffeomorphism between a neighbourhood of the zero section of
$N_{X'_0/\Mloc}$ and a neighbourhood of $X'_0$ in~$\Mloc$ of the form
$U_\eps=\cup_{x\in X'_0}B_x(\eps r)$ for some $\eps>0$.
Denote by $\tilde\phi\in\Omega^3(U_\eps)$ the pull-back of
$\phi_{CY}+R_\theta^* d\eta$ via this diffeomorphism.
Then $F$, as a function of~$v$, becomes equivalent to
$$
v\in W^p_{k}N_{X'_0/\Mloc}\to v^*\tilde\phi\in W^p_{k-1}\Omega^3(X'_0),
$$
for $v$ with sufficiently small in the $W^p_{k}$ norm. Use a uniformizing
coordinate neighbourhood of $w_0$ in $X'_0$ and a finite open cover of the
compact complement by coordinate neighbourhoods and the respective local
trivializations of $N_{X'_0/\Mloc}$ to obtain local expressions for~$v$.
The pull-back $v^*\tilde\phi$ of a given smooth 3-form is then expressed a
cubic polynomial in the partial derivatives of~$v$ with coefficients
smoothly depending on~$v$.
The reader now should have no difficulty to check the smoothness of~$F$
using standard results on Sobolev spaces.
        \end{pf}
The first partial derivative $D_vF(0,0,\theta)$ is a composition of the
bundle isometry $v\mapsto v\corner\phiCY$ and the restriction of
the exterior derivative operator~\eqref{lin.orbi} to a complement of
its kernel. By Theorem~\ref{noobstr}, $D_vF(0,0,\theta)$ maps 
$\tilde{E_0}\oplus W^p_k N_{X'_0/\Mloc}$ isomorphically onto the image
of~$F$. Therefore, the implicit function theorem in Banach spaces applies
to~$F$ and gives for each $\theta_0\in S^1$ a unique smooth family of fields
$v(\eta,\theta)$ in $\tilde{E_0}\oplus W^p_k N_{X'_0/\Mloc}$, for
small $\|d\eta\|$ and $|\theta-\theta_0|$, satisfying
$v(0,\theta)=0$ and $F(\psi(\eta,\theta),\eta,\theta)=0$.
By the compactness of $S^1$ and the uniqueness of the local solutions
$v(\eta,\theta)$, a finite number of these solutions can be patched
together to define a vector field
$(dR_\theta)_0 v(\eta,\theta)$ on $X'_0\times S^1$, for all
$\|d\eta\|<\eps$, where $\eps$ is the smallest of the finitely many
upper bounds coming from the local constructions.

This completes the proof of clauses (i) and (ii) of the theorem and it
remains to establish the regularity clause~(iii).

{\it The interior regularity}. We show that the vector fields
$v(\eta,\theta)$ defined above are smooth on~$X'_0$. In general, a
$(\phiCY+d\eta)$-coassociative submanifold need not be calibrated, so
it is not quite sufficient to fall back on the traditional argument for
minimal submanifolds.
If $v=v_0+v_1\in\tilde{E}_0+W^p_{k} N_{X'_0/\Mloc}$ is a
solution of $F(v,\eta,\theta)=0$ for some $\eta,\theta$
then $v_0$ is smooth and the equation satisfied by $v_1$ can be
written in the form
        \begin{equation}\label{thePDE}
F_0+Gv_1+Q(v_1)=0.
        \end{equation}
Here $Gv_1=d(v_1\lrcorner(\phi_{CY}+d\eta))$ and
$F_0=R_\theta^*d\eta|_{X'_0}+Gv_0$ is a smooth 3-form on $X'_0$.
It follows from the proof of Proposition~\ref{smooth} that the remainder
$Q$ is a cubic polynomial in the first derivatives of $v_1$ with
coefficients smoothly depending on~$v_1$. To deal with this non-linearity
note first that $v_1\in W^p_k N_{X'_0/\Mloc}$, with $k-4/p>1$, thus in
$C^{1,\alpha}$, $\alpha>0$. Then the first factor in each term
$(\nabla v_1)^{\otimes n}\otimes \nabla v_1$ can be considered
as a first order linear differential operator with H\"older continuous
coefficients acting on the second term.
If $\|v_1\|_{C^1}$ is sufficiently small
then the sum of the latter linear operator and $G$ is again an
overdetermined-elliptic operator, so we still have the standard local
regularity estimates~\cite{AGM}. As $v_1$ is already in
$C^{1,\alpha}$ we can apply the usual bootstrapping to show that $v_1$, and
hence also $v$, is smooth.  Then $(dR_\theta)_0\, (v(\eta,\theta))$ defines
a $C^\infty$ vector field on $X'_0\times S^1$.
        \begin{remark}
If the 4-form $*_{\phiCY+d\eta}(\phiCY+d\eta)$ is closed then
$v(\eta,\theta)$ defines a calibrated submanifold and is real analytic.
This is a consequence of a general property of minimal submanifolds
of a Ricci-flat (hence real analytic~\cite{dTK}) Riemannian manifold.
        \end{remark}

{\it $C^{1,\alpha}$-regularity at~$X_0\times S^1$.}
The normal bundle $N_{X'_0/\Mloc}$ is trivial 
and $W^p_k$ sections are $C^{1,\alpha}$ and vanish at~$w_0$ together with first
derivatives when $k-4/p>1$. Multiplying by a smooth cut-off function which is
equal to 1 on $X_0\times S^1$ we obtain an extension of a $L^p_k$ field of
normal vectors to a vector field supported on a
compact neighbourhood of $X_0\times S^1$ on~$\Mloc$. It is clear that such an
extension is smooth except at points of $X_0\times S^1$ where it
is~$C^{1,\alpha}$. It is easy to see that a $W^p_{k}N_{X'_0/\Mloc}$ section
$v(\eta,\theta)$ has an extension to~$\Mloc$ with the same properties as an
$E_0$ component is, by construction, the restriction to~$X_0\times S^1$ of
some smooth vector field on~$\Mloc$.

If $v(\eta)$ has a sufficiently
small $L^p_k$-norm on $X'_0\times\mathrm{pt}$, then we can achieve estimates
$\|v(\eta)\|_{L^p_k(\Mloc)}\|<
\const\|v(\eta,\cdot)\|_{L^p_k(X'_0\times S^1)}$, with constants
independent of~$v$. The estimates 
$\|v(\eta,\theta)\|_{p,k}<\const\|d\eta\|_{p,k}$ follow by differentiating
the identity $F(\psi(\eta,\theta),\eta,\theta)=0$ in $\eta,\theta$ and
taking account of the symmetry of~$F$ in $\theta\in S^1$.
If $\|d\eta\|_{C^1}$ is small then
$\exp_{v(\eta)}$ is a well-defined diffeomorphism
of~$\Mloc$ satisfying the assertions of Theorem~\ref{extn}.

\section{A neighbourhood of the singular fibres}
\label{thrm.B}

In this section we prove Theorem~B and thus complete the proof of the Main
Theorem.
        \begin{varthrm}[Theorem~B]
Let $M$ be a compact 7-manifold with a smooth one-parameter family of
\mbox{$G_2$-structures} given by closed 3-forms $\phi_T\in\Omega^3_+(M)$,
$T>T_1$, defined by the generalized connected sum construction
in~\S\ref{pencils}. Let $\tau_{T}:M\to S^3$ for $T>T_1$ be a coassociative
K3 fibration map defined in Theorem~\ref{approx}.

Then there exists $\eps>0$ so that if $T>T_1$ and an
exact form $d\eta\in\Omega^3(M)$ vanishes on every singular fibre
of~$\tau_{T}$ and $\|d\eta\|_{C^0(M)}<\eps$, relative to the metric
$g(\phi_T)$, then there is a unique smooth vector field
$\tilde{v}_{T}(\eta)$ on~$M$ such that:

(i) $\tilde{v}_{T}$ vanishes on the singular fibres of $\tau_{T}$
and is point-wise orthogonal to each smooth fibre $X$ of $\tau_{T}$
and $\tilde{v}_{T}\corner\phi_T|_X$ is $L^2$-orthogonal to the harmonic
self-dual forms on~$X$ relative to the metric $g(\phi_T)|_X$;

(ii) $\tilde{v}_{T}(\eta)$ depends smoothly on $T$ and $d\eta$ and
$\|\tilde{v}_{T}\|_{C^1}=O(\|d\eta\|_{C^1})$;

(iii) $\phi_T+d\eta$ vanishes on the fibres of
$\tau_{T}\circ\exp_{\tilde{v}_T(\eta)}^{-1}$.
        \end{varthrm}
The main technical issue in the proof of Theorem~B is to establish certain
uniform estimates for families of smooth fibres of~$\tau_T$ 
in a neighbourhood of the singular fibres.
We use the notation of~\S\ref{linear}: $\Mloc=\Wloc\times S^1\subset M$,
$\phiCY\in\Omega^3_+(\Mloc)$, $X_0\subset\Wloc$, $X'_0=X_0\minus\{w_0\}$,
and `holomorphic Morse coordinates' $(z_1,z_2,z_3)$ near~$w_0$.
For $a\in\CX$, let $\xa\subset\Wloc$ denote the fibre of $\tau_T$
which is expressed near $w_0$ by the equation $z_1^2+z_2^2+z_3^2=a^2$
Let $X_{a^2,\theta}=\xa\times\{\theta\}\subset\Wloc\times S^1$
and, as before, identify $\xa=X_{a^2,0}$, for a marked point $0\in S^1$.
The fibre $X_{a^2,\theta}$ is well-defined for any small $|a|$, say
$|a|<\rho$, and is non-singular for $a\neq 0$.
Throughout this section any $X_{a^2,\theta}$ is considered with the
metric $\ga $ induced by restriction of $g(\phiCY)$ unless stated otherwise.
For each compact smooth coassociative submanifold~$X\subset M$, define
a decomposition
        \begin{gather*}
L^p_k(N_{X/M})=
\Gamma^h(N_{X/M})\oplus\Gamma^\bot(N_{X/M})_{p,k},
\intertext{where}
\Gamma^h(N_{X/M})=\{v\in\Gamma(N_{X/M})\,|\; 
(v\lrcorner\phi_T)\in\HH^+(X)\},
\\
\Gamma^\bot(N_{X/M})=\{v\in\Gamma(N_{X/M})\,|\; 
(v\lrcorner\phi_T)\bot_{L^2}\HH^+(X)\},
        \end{gather*}
and $\Gamma^\bot(N_{X/M})_{p,k}$ means the completion of
$\Gamma^\bot(N_{X/M})$ in the $L^p_k$ norm (as before, we require
$1<k-4/p<2$). Denote
$$
\eus{Y} = \{\eta\in \Omega^2(\Mloc)\,|\; d\eta|_{X_{0,\theta}}=0
\text{ for each }\theta\in S^1\}.
$$
The proof of Theorem~B builds up on the following.
        \begin{prop}\label{near.orbi}
There is an $\eps_0>0$ so that whenever $\eta\in\eus{Y}$ and
$\|d\eta\|_{C^1(\Mloc)}<\eps_0$ the following holds.
For each $0<|a|<\rho$ and $\theta\in S^1$ there is a unique smooth
$v=v_{a^2,\theta}(\eta)\in\Gamma^\bot(N_{X_{a^2,\theta}/\Mloc})$ such
that $v_{a^2,\theta}(0)=0$, $v_{a^2,\theta}(\eta)$ depends smoothly on
$\eta$, and $\phi_{CY}+d\eta$ vanishes on the submanifold
$\exp_{v_{a^2,\theta}(\eta)}(X_{a^2,\theta})$.

For a fixed $\eta\in\eus{Y}$ with $\|d\eta\|_{C^1(\Mloc)}<\eps_0$, the family
$v_{a^2,\theta}(\eta)$, for $0<|a|<\rho$, $\theta\in S^1$ extended by zero
over $X_{0,\theta}$,
defines a vector field $v(\eta)$ which is $C^1$ on a neighbourhood
$\Mloc$ of~$X_0\times S^1$ and smooth on $\Mloc\setminus X_0\times S^1$.
Furthermore, $\|v(\eta)\|_{C^1(\Mloc)}<K\|d\eta\|_{C^1(\Mloc)}$, for some
constant~$K$ independent of $\eta$.
        \end{prop}
        \begin{pf}
For each $a\neq 0$, $\theta\in S^1$ we want to obtain
$v_{a^2,\theta}(\eta)$ as $v(0,\eta)$, where $v(\beta,\eta)$ is the
solution to implicit function problem $F(v(\beta,\eta),\beta,\eta)=0$ for
the extended McLean's map for $\xat$
        \begin{multline}\label{extended}
F:\; (v,\beta,\eta)
\in\Gamma^\bot(N_{X_{a^2,\theta}/\Mloc})_{p,k}
\oplus\Gamma^h(N_{X_{a^2,\theta}/\Mloc})\oplus\eus{Y}
\to\\
\exp^*_{v+\beta}(\phiCY+d\eta)|_{X_{a^2,\theta}}
\in L^p_{k-1}\Omega^3(X_{a^2,\theta})\cap\Ker d,
        \end{multline}
where $k-4/p>1$ we use the completion in the $C^{3}$ norm for $\eta\in\eus{Y}$.
That $F$ is a smooth map between the indicated Banach spaces follows by
the argument of Proposition~\ref{smooth}, up to a change of notation.
We have $F(0,0,0)=0$ and
the derivative $(D_vF)_{0,0,0}v=d(v\corner\phiCY)|_{X_{a^2,\theta}}$ is an
isomorphism of Banach spaces by standard Hodge theory as
$H^3(X_{a^2,\theta},\RE)=0$ for a K3 surface $X_{a^2,\theta}$.

The implicit function theorem in Banach spaces gives a unique family of
sections $v(\beta,\eta)$, with $v(0,0)=0$, defined for
$\max\{|\beta|,\|d\eta\|_{C^1(X_{a^2,\theta})}\}<\delta$ say,
but this $\delta$ may in general depend on~$a$. (We can choose $\delta$ 
independent of~$\theta$ by the $S^1$-symmetry of~$\phiCY$.)
To keep track of the relation between constants appearing in
the estimates we use.
        \begin{prop}[Implicit function theorem in Banach spaces]\label{IFT}
Suppose that a smooth map \mbox{$f:E=E_1\oplus E_2\to F$} between Banach
spaces has an expansion
$$
f(\xi_1,\xi_2)=(D_1f(0,0))\xi_1+(D_2f(0,0))\xi_2+Q(\xi_1,\xi_2),
$$
so that $A=D_1f(0,0):E_1\to F$ is an isomorphism of Banach spaces and for
$\xi,\zeta\in E$,
$$
\|A^{-1}Q(\xi)-A^{-1}Q(\zeta)\|<C(\|\xi\|+\|\zeta\|)\|\xi-\zeta\|,
$$
for some constant $C$. Then there
exists a uniquely determined smooth function
$\phi:B_{2,\delta}\to B_{1,\delta}$, $\phi(0)=0$, where
$B_{i,\delta}=\{\xi_i\in E_i: \|\xi_i\|<\delta\}$ with $\delta=(4C)^{-1}$,
so that all zeros of $f$ in $B_{1,\delta}\times B_{2,\delta}$ are of the
form $(\phi(\xi_2),\xi_2)$.
        \end{prop}
The proof of Proposition~\ref{IFT} is a standard application of the
contraction mapping principle.

We require a lower bound on
the linearization $(D_1 F)_0$ in $v$ and an upper bound on the quadratic
remainder for~$F$ so that the constant $(4C)^{-1}$ in Proposition~\ref{IFT}
for the map $F$ on a manifold $\xa$ is greater than the norm of
$d\eta|_{X_{a^2,\theta}}$ for any sufficiently small $|a|\neq 0$.
As noted above, the estimates on $X_{a^2,\theta}$, can be taken
independent of $\theta\in S^1$.

We begin with the linear part.
        \begin{prop}\label{lower.bound}
There exists a $\rho>0$ such that for any
$v\in\Gamma^\bot(N_{X_{a^2,\theta}/\Mloc})_{p,k}$
with $0<|a|<\rho$,
$$
\|v\|_{p,k}<C_{p,k}\|d(v\corner\phiCY)|_{\xat}\|_{p,k-1}
$$
with a constant $C_{p,k}$ independent of $v$ or~$a$.
        \end{prop}
        \begin{pf}
Since the linear map and the norms are symmetric in $\theta$ we drop
$\theta$ from the notation.
Consider a sphere $|z|=|a|^{1/2}$ in the coordinate neighbourhood of $w_0$
in~$\Wloc$, where $|z|^2=\sum_{j=1}^3|z_j|^2$. If $|a|$ is sufficiently
small then this sphere intersects~$\xa$ and we can write
$\xa=\xa^-\cup\xa^+$, where
$\xa^-=\xa\minus\{|z|\le|a|^{1/2}\}$ and
$\xa^+=\xa\cap\{|z|<2|a|^{1/2}\}$ are two open submanifolds of~$\xa$.
As $a\to 0$ the metric on $\xa^-$ is asymptotic in $C^\infty$ to the
metric on $X_0\minus\{|z|\le|a|^{1/2}\}$. The other piece $\xa^+$ with
metric rescaled by a constant factor $|a|$ is asymptotic in~$C^\infty$
to a complex surface $\Sigma=\{\sum_{j=1}^3 z_j^2=1,\;\; |z|<|a|\}$ in
$\CX^3$ with the metric on~$\Sigma$ induced by a Hermitian inner product
on $T_{w_0}\Wloc\cong\CX^3$ defined by the K\"ahler form~$\omega(w_0)$.

The decomposition $\xa=\xa^-\cup \xa^+$
can be thought of as a generalized connected sum $X'_0\#_{\RP^3}\Sigma$ of
two manifolds taken at their ends $\RE_{>0}\times\RP^3$.
The metric $\ga$ on the connected sum $X_a$ is $C^\infty$ asymptotic to a
metric smoothly interpolating between metrics on compact subsets of $X'_0$
and~$\Sigma$.

A family of Riemannian 4-manifolds $(\xa,\ga)$ is an instance of the gluing
construction studied in~\cite[\S 2]{KS}. The estimate that we are
interested in is equivalent to a lower bound on the operator
        \begin{equation}\label{d}
d:L^p_k\Omega^+(\xa)\to L^p_{k-1}\Omega^3(\xa).
        \end{equation}
The `main estimate' proved in~\cite[\S~4.1]{KS} deals with the invertibility
of an elliptic differential operator on a generalized connected sum, as
above, when the coefficients of this operator are obtained by gluing the
coefficients of Fredholm elliptic operators on the two manifolds with
cylindrical ends. The result, in particular, asserts a uniform lower bound
on a subspace of finite codimension in the $L^p_k$ domain on a connected
sum. Inspection of the proof of this uniform lower bound shows that it
remains valid for the overdetermined elliptic operator~\eqref{d} on
$(\xa,\ga)$ provided that the $L^2$ kernels of the respective exterior
derivatives on self-dual forms on $X'_0$ and $\Sigma$ are
finite-dimensional. Then the operator~\eqref{d} admits a lower bound,
independent of~$a$ as $a\to 0$, on the complement of a finite-dimensional
subspace in $L^p_k\Omega^+(\xa)$. The dimension of this space is the sum of 
dimensions of the spaces of closed self-dual $L^2$ forms on $X'_0$
and~$\Sigma$.

These latter spaces are contained in the $L^2$ kernels of the Laplacian
on self-dual forms on $X'_0$ and $\Sigma$. The dimensions can be computed
in the asymptotically cylindrical metrics in the conformal class of
$X'_0$ and $\Sigma$ as the $L^2$ norm of 2-forms is conformally invariant
in dimension~4. Recall from Proposition~\ref{1-1} that the $L^2$ kernel is 
3-dimensional for $\Omega^+(X'_0)$. The 4-manifold $\Sigma$ is
diffeomorphic to the total space of the tangent bundle $TS^2$, so
$b^+(\Sigma)=0$ and the same argument as in Proposition~\ref{1-1} shows
that the Laplacian on $L^p_k\Omega^+(\Sigma)$ is injective. Thus the uniform
lower bound for~\eqref{d} holds on the complement of a 3-dimensional subspace
in~$\Omega^+(\xa)$. It is {\it a posteriori} clear that this 3-dimensional
space can be taken to be~$\HH^+(\xa)$.
        \end{pf}
The argument of Proposition~\ref{smooth} shows that the non-linear
remainder $Q$ of the map~$F$ consists of polynomial terms in
$\nabla v$ with coefficients smooth in~$v$. This easily gives an
upper bound on~$Q$ independent of $a,\theta$ for small~$|a|$.
Therefore, by Propositions~\ref{IFT} and~\ref{lower.bound} there is an
$\eps>0$ independent of $a,\theta$ so that if $\eta\in\Omega^2(\Mloc)$
satisfies $\|d\eta|_{\xa,\theta}\|_{C^1(X_{a^2,\theta})}<\eps$ then
the implicit function theorem defines sections $v(0,\eta)$
of $N_{\xat/\Mloc}$ for every small~$|a|$.

The metric on the submanifold $X_{a^2,\theta}$ is independent of $\theta\in
S^1$ and in the next result we once again drop $\theta$ from the notation. 
It will be convenient to use H\"older norms commensurable with the $L^p_k$
norms that we required.
        \begin{prop}\label{et}
For each $\xa$ with $a\neq 0$, we have
$$
\|d\eta|_{\xa}\|_{C^{1,\alpha}(\xa)}<
\lambda_{k,\alpha}\|d\eta\|_{C^{1,\alpha}(\Mloc)}|a|^{1-\alpha},
\quad
0\le\alpha<1,
$$
with a constant $\lambda_{k,\alpha}$ independent of~$\eta$ or~$a$.
        \end{prop}
Estimates similar to that in Proposition~\ref{et} are proved in
\cite{joyce-cone-3} and~\cite{lotay-gafa}. In the present situation, we
have an explicit algebraic local model for the submanifolds $\xa$ near the
singular point of $X_0$ and some details are simplified.
        \begin{pf}
Recall from the proof of Proposition~\ref{lower.bound} 
the decomposition $\xa=\xa^-\cup \xa^+$ as a generalized connected sum
and denote by $(d\eta)^{\pm}$ the restrictions of $d\eta$ to the respective
pieces $\xa^\pm\subset\xa$.
Recall also that $\xa^-$ is diffeomorphic to $X_0\minus\{|z|\le|a|^{1/2}\}$.
As $X'_{a^2}=X_0\minus\{|z|\le|a|^{1/2}\}$ is compact and non-singular, its
injectivity radius of is bounded away from zero, furthermore, it is not
difficult to check that the injectivity radius is $O(|a|$, as $|a|\to 0$.
We deduce that for $|a|$ sufficiently small, each $\xa^-$ is a well-defined
deformation of $X'_{a^2}$ defined by a section $\nu_{a^2}$ of the normal
bundle of the latter submanifold, using the exponential map for
$g(\phiCY)$. The $\nu_{a^2}\to 0$ depends smoothly on $a$ and vanishes when
$a=0$. As $d\eta$ is smooth and vanishes on $X'_0$ and the metric on $\xa^-$
converges in $C^\infty$ to the metric on $X'_0$, we find that the $C^k$
norm of $(d\eta)^-$ is a smooth function of~$a$. This function vanishes
at $a=0$, hence is bounded by a constant multiple of
$\|d\eta\|_{C^k(\Mloc)}|a|$.

The form $(d\eta)^+$ is defined on a coordinate neighbourhood
$\{|z|<|a|^{1/2}\}\times S^1$ of $w_0$ in $\Mloc$. We may assume, by
rescaling $z\in\CX^3$ if necessary, that the matrix of the inner product
defined by $g(\phiCY)$ on $T_{w_0}\Wloc$ has determinant~1. The metric
$g(\phiCY)$ on the coordinate neighbourhood may be written as $g_0+O(|z|)$,
where $g_0$ denotes the value of $g(\phiCY)$ at~$w_0$ and $O(|z|)$ is
estimated independently of~$\theta$. A similar local expansion is valid for
the derivatives of $g(\phiCY)$.

The parameterizations of local submanifolds
$\xa^+$ may be taken to be homogeneous of order 1 in~$a$. The
coefficients of the $k$-th derivatives of induced metric $g(\xa^+)$ on
$\xa^+$ then have a local extension with the leading 
term homogeneous of order $2-k$ in $a$ and a remainder $O(|a|^{3-k})$ as
$a\to 0$ ($d\theta$ vanishes on each $\xa$). The difference between the
restriction of the Levi--Civita connection of $g(\phiCY)$ on the ambient
$\Mloc$ to $\xa^+$ and the Levi--Civita connection of the metric $g(\xa^+)$
induced by $g(\phiCY)$ on $\xa$ is determined by second fundamental form of
$g(\xa^+)$. Recall that for a submanifold defined by submersion $\tau_T$
the second fundamental form is the quotient of the Hessian of~$\tau_T$ and
the gradient of~$\tau_T$. Thus its point-wise norm has leading term
homogeneous of order $a^{-1}$, as measured by $g(\xa^+)$. As $d\eta$
vanishes on $X_0$, this 2-form vanishes on $T_{w_0}\Wloc$, the Zariski
tangent space of $X_0$ at $w_0$, so the coefficients of $d\eta|_{\Wloc}$
have zeros of order two at~$w_0$ and the covariant derivative has zero of
order one. It follows that for $0<|a|<\rho$ with sufficiently
small~$\rho$, the $C^{1,\alpha}$ norm of $d\eta|_{\xa^+}$ measured using
the Levi--Civita of $g(\xa^+)$ is $O(\|d\eta\|_{C^k(\Mloc)}|a|^{1-\alpha})$.

As $\|d\eta|_{\xa}\|_{C^k(\xa)}\le
\max\{\|d\eta^-|_{\xa^-}\|_{C^k(\xa^-)},\|d\eta^+|_{\xa^+}\|_{C^k(\xa^+)}\}$
the proposition is proved.
        \end{pf}
The local deformation $\exp_{v(\beta,\eta)}(X_{a^2,\theta})$ is
well-defined if $v(\beta,\eta)$ is small in the uniform norm so that the
tubular neighbourhood theorem is applicable to~$X_{a^2,\theta}$. It is not
difficult to see that, for any small $a$, a suitable neighbourhood of the
zero section of $N_{\xat/\Mloc}$ is $\|v\|_{C^0}<\eps_1\,|a|$, 
for a sufficiently small constant $\eps_1>0$ independent of~$a$.

For a small $d\eta$ we can estimate on $\xat$, \
$\|v(0,\eta)\| = O(\|(D_2v)_{0,0}\|\,\|d\eta|_{\xat}\|)$.
As $(D_2v)_{0,0}=(D_1F)_{0,0,0}^{-1}(D_3F)_{0,0,0}$
this gives
        \begin{equation}\label{est}
\|v(0,\eta)\|_{C^1(X_{a^2,\theta})} = O(|a|),\qquad a\to 0,
        \end{equation}
using Proposition~\ref{lower.bound} and Lemma~\ref{et}.
In particular, for a sufficiently small $d\eta$ the deformations 
$\exp_{v(0,\eta)}(X_{a^2,\theta})$ are well-defined for every small~$|a|$.

The regularity results for the coassociative submanifolds~\cite{calibrated}
imply that $v(\beta,\eta)$ is a smooth section over $\xat$, as $v(\beta,\eta)$
is already in~$C^1$.  By construction, $\exp_{v(\beta,\eta)}(\xat)$ coincides
with $\exp_{v(0,\eta)}(\exp_{v(\beta,0)}(\xat))$. Therefore, the local
families $v(\beta,\eta)$ can be patched together to define a smooth
vector field $v(\eta)$ on $(\cup_{0<|a|<\rho}\xa)\times S^1$, for some $\rho>0$.
This vector field has a Lipschitz
continuous extension by zero over~$\cup_{\theta\in S^1}X_{0,\theta}$.

Differentiating the identity $F(v(\beta,\eta),\beta,\eta)=0$ in~$\beta$
at $\beta=0$ we find that the derivative of $v(\eta)$ in the directions
orthogonal to the fibres $\xat$ satisfies
$$
(D_1F)_{v(0,\eta),0,\eta}(D_1v)_{0,\eta}+(D_2F)_{v(0,\eta),0,\eta}=0.
$$
Using explicit expressions for the 
derivatives of~$F$ and the vanishing of $\phi_{CY}$ on $\xat$ we obtain
$$
d\bigl(((D_1v)_{0,\eta}\beta)\corner \phiCY\bigr)\bigl|_{\xat}=
-d(\beta\corner d\eta)\bigr|_{X_{a^2,\theta}}.
$$
The operator $d(\cdot\corner\phiCY)|_{X_{a^2,\theta}}$ is bounded below
independent of $a$ by Proposition~\ref{lower.bound} and we deduce
that $\|(D_1v)_{0,\eta}\|_{C^0(\xat)}< K\|d\eta\|_{C^1(\Mloc)}$, $a\neq 0$.

We adapt the same method to show that 
$v(\eta)$ is continuously differentiable at any point in $X_0\times S^1$.
Differentiating the identity $F(v(\beta,\eta),\beta,\eta)=0$ in~$\beta$
at $\beta=0$ we find that the derivative of $v(\eta)$ in the directions
orthogonal to the fibres $\xat$ satisfies
$$
(D_1F)_{v(0,\eta),0,\eta}(D_1v)_{0,\eta}+(D_2F)_{v(0,\eta),0,\eta}=0.
$$
The self-dual harmonic forms on $\xat$
are spanned by the K\"ahler form $\omega|_{\xat}$ and the real and imaginary
parts of the $(2,0)$-form $(\p\tau)^\sharp\corner\Omega$ induced by the
Calabi--Yau structure $(\omega,\Omega)$ on $\Wloc\times\{\theta\}$.
The definition of $\Gamma^h(N_{\xat/M})$ therefore can be extended to
the smooth subset $X'_{0,\theta}\subset X_{0,\theta}$. The resulting
$\Gamma^h(N_{X'_{0,\theta}/M})$ can be thought of as a limit of
$\Gamma^h(N_{\xat/M})$ as $a\to 0$, in the sense that the 4-manifold
$X'_{0,\theta}$ is diffeomorphic to an open subset of $\xat$ and the metric
on $X'_{0,\theta}$ is a $C\infty$ limit of the metrics on this subset.
Using explicit expressions for the 
derivatives of~$F$ and taking the limit as $a\to 0$ we obtain
$$
d\bigl(((D_1v)_{0,\eta}\beta)\corner \phiCY\bigr)\bigl|_{X'_{0,\theta}}=
-d(\beta\corner d\eta)\bigr|_{X'_{0,\theta}}
$$
The operator $d(\cdot\corner\phiCY)|_{X'_{0,\theta}}$ is injective
on $\Gamma^\bot(N_{X'_{0,\theta}/\Mloc})$. Applying a left
inverse we determine $(D_1v)_{0,\eta}$ at any point in $X'_{0,\theta}$.
This shows the continuity of $(D_1v)_{0,\eta}$ away from critical point
of~$\tau_T$. The estimate~\eqref{est} shows that the derivative of
$v(\eta)$ is Lipschitz continuous at each point in $w_0\times S^1$, thus
$v(\eta)$ is $C^1$-regular on all of~$\Mloc$. The higher order derivatives
of $v(\eta)$ can be handled by a similar method but the  expressions become
cumbersome; we omit the details.

We estimated the `intrinsic' first derivatives of $v(\eta)$ on $\xat$ and
in the transverse directions by $\|d\eta\|_{C^1(\Mloc)}$. To obtain the
desired estimate of the $C^1$ norm of $v(\eta)$ on~$\Mloc$ we recall that
the second fundamental form of~$\xat$ is $O(|a|^{-1})$, as $a\to 0$, near the
critical points of $\tau_T$ and bounded on the complement of a neighbourhood
of the critical points. It follows that 
$\|v\|_{C^1(\Mloc)}=O(\|d\eta\|_{C^1(\Mloc)})$.
This completes the proof of Proposition~\ref{near.orbi}.
        \end{pf}

The singular fibres of~$\tau_T$ occur in finitely many $S^1$-families.
Applying Proposition~\ref{near.orbi} to each $S^1$-family we obtain for any
$\eta\in\Omega^2(\Mloc)$ with small $\|d\eta\|_{C^1}$ a vector field
$\tilde{v}_T(\eta)$ satisfying the assertions of Theorem~B except that
$\tilde{v}_T(\eta)$ is only defined on a neighbourhood $U$ of the
singular fibres of~$\tau_T$.

The implicit function argument for the extended McLean's map
$F(v,\beta,\eta)$ as in~\eqref{extended} for a fibre $X$ of $\tau_T$
in $M\!\minus U$ gives a uniquely determined smooth local family
$v(\beta,\eta)$ of sections in $\Gamma^\bot(N_{X/M})$, so that
$v(0,0)=0$ and $\phi_T+d\eta$ vanishes on $\exp_{v(\beta,\eta)}(X)$.

The estimates on $v(\beta,\eta)$ are similar to those in
Proposition~\ref{near.orbi} but easier as the curvature of the fibres in
$M\!\minus U$ is bounded independent of the fibre.
The complement $M\!\minus U$ is a family of non-singular fibres
of~$\tau_T$. Recall from \S\ref{pencils} that the metric on each of
these fibres is up to a small deformation a metric in a {\em compact}
closure of the finite-dimensional family of metrics on the fibres of
holomorphic fibrations $\tau^{(j)}$, $j=1,2$, restricted to the
asymptotically cylindrical ends of~$W_j$. The deformation is bounded,
independent of $T$ or the fibre, in $C^k$ norm for each~$k$.
It follows that the lower bound on $(D_1F)_{0,0,0}$ and the upper bound on
the quadratic remainder can be taken independent of $X\subset M\!\minus U$.
Then the sections $v(\beta,\eta)$ are defined for
$\|d\eta\|_{C^1(M)}<\eps$ where $\eps$ can be taken independent of $T$ or
the fibre $X$. 

Using once again the property that for a fixed $\beta$ the sections
$v(\beta,\eta)$ define deformations of $\exp_{v(\beta,0)}(X)$, a fibre
near $X$, we obtain that $v(\beta,\eta)$ induces a vector field satisfying
the assertions of Theorem~$B$ except that it is only defined on an open
neighbourhood of $X$ in~$M$.  The vector fields obtained from
$v(\beta,\eta)$ for different $X$ agree on the overlaps of their domains
by the uniqueness part of the implicit function theorem.
As $M\!\minus U$ is compact an extension of $v(\eta)$ from $U$ to $M$ is
obtained by patching with a finite number of these local vector fields.

\end{document}